\documentclass[11pt,a4paper]{article}
\usepackage[numbers,sort&compress]{natbib}
\usepackage[colorlinks,
            linkcolor=blue,
            anchorcolor=green,
            citecolor=magenta
            ]{hyperref}
\usepackage{mathrsfs,amssymb,amsfonts}
\usepackage{graphicx}
\usepackage{amsmath,amssymb,latexsym,color}
\usepackage[mathscr]{eucal}
\usepackage[numbers]{natbib}
\usepackage{bm}
\usepackage[top=1in, bottom=1in, left=1.25in, right=1.25in]{geometry}
\makeatletter

\newcommand{\Rmnum}[1]{\expandafter\@slowromancap\romannumeral #1@}
\makeatother
\newtheorem{thm}{Theorem}[section]
\newtheorem{defin}[thm]{Definition}

\newtheorem{lem}[thm]{Lemma}

\newtheorem{example}{Example}
\newtheorem{remark}{Remark}[section]

\newcommand{\qed}{\hfill\Box\medskip}
\usepackage{CJK}
\usepackage{geometry}
\begin{document}
\begin{CJK*}{GBK}{song}
 \setlength{\baselineskip}{14pt}
\renewcommand{\abovewithdelims}[2]{
\genfrac{[}{]}{0pt}{}{#1}{#2}}

\title{\bf Fault-free Hamiltonian cycles in balanced hypercubes with conditional edge faults }

\author{Pingshan Li \quad  Min Xu\footnote{Corresponding author. \newline {\em E-mail address:} xum@bnu.edu.cn (M. Xu).}\\
{\footnotesize   \em  Sch. Math. Sci. {\rm \&} Lab. Math. Com. Sys., Beijing Normal University, Beijing, 100875,  China} }
 \date{}
 \date{}
 \maketitle

\begin{abstract}
The balanced hypercube, $BH_n$, is a variant of hypercube $Q_n$.
Zhou et al. [Inform. Sci. 300 (2015) 20-27]  proposed an interesting  problem that whether there is a fault-free  Hamiltonian cycle in $BH_n$ with each vertex incident to at least two fault-free edges. In this paper, we consider this problem and show that each fault-free edge lies on a fault-free Hamiltonian cycle in $BH_n$ after no more than $4n-5$ faulty edges occur if each vertex is incident with at least two fault-free edges for all $n\ge 2$. Our result is optimal
 with respect to the maximum number of tolerated edge faults.

\medskip
\noindent {\em Key words:} Balanced hypercubes; Hypercubes; Hamiltonian cycles; Fault-tolerance.

\medskip
\end{abstract}

\section{Introduction}
In the field of parallel and distributed systems, interconnection networks are an important research area. Typically, the topology
of a network can be represented as a graph in which the vertices represent processors and the edges represent communication links.

For graph definitions and notations, we follow \cite{J.A.Bondy}.
A graph $G$ consists of a vertex set $V(G)$ and an edge set $E(G)$, where an edge is an unordered pair of distinct vertices of $G$. A graph $G$ is called bipartite if its vertex set can  be partitioned into two parts $V_1, V_2$  such that every edge has one endpoint in $V_1$ and one in $V_2$. A vertex $v$ is a neighbor of $u$ if $(u, v)$ is an edge of $G$, and $N_{G}(u)$ denotes all the neighbors of $u$ in $G$.
A path  $P$ of length $\ell$ from $x$ to $y$ is a finite sequence of distinct vertices $\langle v_0, v_1, \cdots, v_{\ell}\rangle$ such that $x=v_0, y=v_{\ell}$, and $(v_i, v_{i+1})\in E$ for $0\le i\le \ell-1$.
We also denote the path $P$ as $\langle v_0, v_1, \cdots, v_i, Q, v_j, v_{j+1}, \cdots, v_{\ell}\rangle$, where $Q$ is the path $\langle v_i, v_{i+1}, \cdots, v_j\rangle$. A cycle $C$ of length $\ell+1$ is a closed path$\langle v_0, v_1, \cdots, v_{\ell}, v_0\rangle$.

A path (resp.,  cycle) is called a Hamiltonian path (resp.,  cycle) if it contains all the vertices of $G$. A graph $G$ is said to be Hamiltonian if there is a Hamiltonian cycle.
A graph $G$ is said to be Hamiltonian connected if there is a Hamiltonian path between any two vertices of $G$. A bipartite graph is Hamiltonian laceable if there is a Hamiltonian path between any two vertices in different bipartite sets.

The Hamiltonian property is one of the major requirements in designing network topologies since a topology structure containing Hamiltonian paths or cycles can efficiently simulate many algorithms that are designed on linear arrays or rings \cite{Wangwenqing}.

It is important to consider fault-tolerance in  networks since faults may occur in  real  networks. Two fault models have been studied in many well-known networks, one is the random fault model, which means that faults may occur anywhere without any restriction, see, for example,  \cite{Xujunming,Mameijie,Yangmingchen3,Sunsainan}.
The other is the conditional fault model, which assumes that the fault distribution is limited. For example, some studies on two or more non-faulty edges incident to each vertex can be found in \cite{Wanghailiang, Yangdawei, Liujiajie,Hunghaoshun,Xujunming2,Hsiehsunyuan,Chenjhengcheng}.

The hypercube network has been proved to be one of the most popular interconnection networks because it possesses many excellent properties such as a recursive structure, regularity, and symmetry.
The balanced hypercube, proposed by Huang and Wu \cite{Huang}, is  a  hypercube variant. Similar to hypercubes, balanced hypercubes are bipartite graphs \cite{Huang} that are vertex-transitive \cite{Huang2} and edge-transitive \cite{Zhou}. Balanced hypercubes are superior to hypercubes in that they have a smaller diameter than hypercubes\cite{Huang2}.

The balanced hypercube, $BH_n$, has been studied by many researchers in recent years.
Xu  et al. \cite{Xumin} proved that $BH_n$ is edge-bipancyclic and Hamiltonian laceable.
Yang \cite{Mingchengyang} proved that  $BH_n$ is bipanconnected.
Yang  \cite{Mingchengyang2} also demonstrated that the super connectivity of $BH_n$ is $4n-4$ and the super edge-connectivity of $BH_n$ is $4n-2$ for $n\ge 2$.
L\"{u} et al. \cite{Huazhonglu} proved that $BH_n$ is hyper-Hamiltonian laceable.
Cheng et al. \cite{Dongqincheng} proved that $BH_n$ is $(n-1)$-vertex-fault-tolerant edge-bipancyclic.
Hao et al. \cite{R.X.Hao} showed that there is a fault-free Hamiltonian path between any two adjacent vertices in $BH_n$ with $2n-2$ faulty edges.
Cheng et al. \cite{Dongqincheng2}  proved that $BH_n$ is $2n-3$ edge-fault-tolerant $6$-edge-bipancyclic  for all $n\ge2$.
Zhou et al. \cite{Qingguozhou} proved that $BH_n$ is $2n-2$ edge-fault-tolerant Hamiltonian laceable, and they proposed an  interesting problem that whether  there is a  fault-free Hamiltonian cycle in $BH_n$ with each vertex incident to at least two fault-free edges. In this paper, we consider this problem and  show that each fault-free edge lies on a fault-free Hamiltonian cycle in $BH_n$ after no more than $4n-5$ faulty edges occur if each vertex is incident with at least two fault-free edges for all $n\ge 2$.  Our result is optimal with respect to the maximum number of tolerated edge faults.

The rest of this paper is organized as follows. In Section $2$, we introduce  two equivalent definitions  of balanced hypercubes  and discuss some of their properties.  In section $3$, we introduce some lemmas used in the  proof of the main result and prove the main result.  Finally, we conclude this paper and give an example to show that our result is optimal in Section $4$ .

\section{Balanced hypercubes}

\begin{figure}[!htbp]
  \centering
  \includegraphics[width=0.6\textwidth]{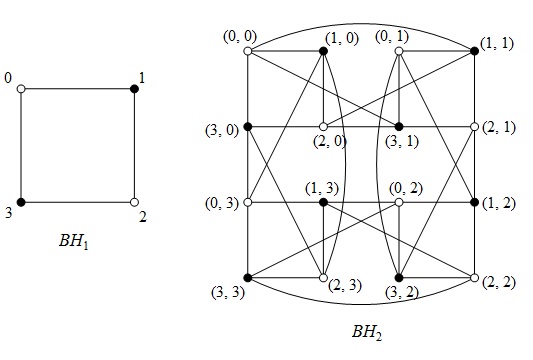}\\
  \caption{Illustration of $BH_1$ and $BH_2$}\label{BH2}
\end{figure}

Wu and Huang  \cite{Huang} present two equivalent definitions of $BH_n$ as follows:

\begin{defin}  An $n$-dimensional balanced hypercube $BH_n$ has $2^{2n}$ vertices, each labeled by an $n$-bit string  $(a_0, a_1, \cdots, a_{n-1})$, where $a_i\in \{0, 1, 2, 3\}$ for all $0\le i\le n-1$. An arbitrary vertex $(a_0, a_1, \cdots, a_{i-1}, a_{i}$, $a_{i+1}, \cdots, a_{n-1})$ is adjacent to the following $2n$ vertices:
\begin{center}
  $\begin{array}{l}
          (1)~~ ((a_0\pm1)\mod\ 4, a_1, \cdots, a_{i-1}, a_{i}, a_{i+1}, \cdots, a_{n-1})~where~1\le i\le n-1, \vspace{0.5ex} \\
           (2)~~ ((a_0\pm1)\mod\ 4, a_1, \cdots,  a_{i-1}, (a_{i}+(-1)^{a_0})\mod\ 4, a_{i+1}, \cdots, a_{n-1})~where~1\le i\le n-1.
         \end{array}
  $
 \end{center}
In $BH_n$, the first coordinate $a_0$ of vertex $(a_0, a_1, \cdots, a_{n-1})$ is called the inner index, and the second coordinate $a_i(1\le i\le n-1)$ is called the $i$-dimension index. From the definition, we have  $N_{BH_n}((a_0, a_1, \cdots, a_{n-1}))=N_{BH_n}((a_0+2, a_1, \cdots, a_{n-1}))$.
Figure \ref{BH2} shows two balanced hypercubes of dimensional one and two.
\end{defin}

 Briefly, we assume that `$+, -$' for the coordinate of a vertex is an operation with mod 4 in the remainder of the paper.
 Let $X_{j, i}=\{(a_0, a_1, \cdots, a_{j-1}, a_{j}, a_{j+1}, \cdots, a_{n-1})\mid a_k\in \{0, 1, 2, 3\}, 0\le k\le n-1, a_j=i\}$ for $1\le j\le n-1$ and $ i\in \{0, 1, 2, 3\}$ and let $BH^{j, i}_{n-1}=BH_n[X_{j, i}]$. Then, $BH_n$ can be divided into four copies: $BH^{j, 0}_{n-1}, BH^{j, 1}_{n-1}, BH^{j, 2}_{n-1}$, and $BH^{j, 3}_{n-1}$ where $BH^{j, i}_{n-1}\cong BH_{n-1}$ for $i=0, 1, 2, 3$ \cite{Dongqincheng}. We use $BH^i_{n-1}$ to denote  $BH^{n-1, i}_{n-1}$ for $i=0,1, 2, 3$ .

\begin{defin} An $n$-dimensional balanced hypercube  $BH_n$ can be constructed recursively as follows:
\begin{enumerate}
  \item $BH_1$ is a cycle of length four with vertex-set $\{0, 1, 2, 3\}$.
  \item $BH_{n}$ is a construct from four copies of $BH_{n-1}: BH^0_{n-1}, BH^1_{n-1}, BH^2_{n-1}$, and  $BH^3_{n-1}$. Each vertex $(a_0, a_1, \cdots$, $ a_{n-2}, i)$ has two extra adjacent vertices:

$\begin{array}{l}
  ~(1)~ In~ BH^{i+1}_{n-1}: (a_0 \pm 1, a_1, \cdots,  a_{n-2}, i+1)~if~a_0~is~ even\vspace{1.0ex}. \\
   ~(2)~ In~ BH^{i-1}_{n-1}: (a_0 \pm 1, a_1, \cdots, a_{n-2}, i-1)~if~a_0~is~ odd.
 \end{array}
$

\end{enumerate}

\end{defin}

Since $BH_n$ is a bipartite graph,  $V(BH_n)$ can be divided into two disjoint parts. Obviously, the vertex-set $V_1=\{a=(a_0, a_1, \cdots, a_{n-1})\mid a\in V(BH_n)$ with $ a_0$  odd$\}$ and $V_2=\{a=(a_0, a_1, \cdots, a_{n-1})\mid a\in V(BH_n)$ with $ a_0$   even$\}$ form the desired partition. We use black nodes to denote the vertices in $V_1$ and white nodes to denote the vertices in $V_2$.

Let $(u, v)$ be an edge of $BH_n$. If $u$ and $v$ differ only with regard to the inner index, then $(u, v)$ is  said to be a 0-dimension edge. If $u$ and $v$ differ not only in terms of the inner index but also with regard to the $i$-dimension index, then $(u, v)$ is called the $i$-dimension edge. We use $\partial D_d(0\le d\le n-1)$ to denote the set of all $d$-dimension edges.

There are some known properties about $BH_n$.

\begin{lem}{\rm (\cite{Huang2, Zhou})}\label{L1} The balanced hypercube $BH_n$ is vertex-transitive and edge-transitive.
\end{lem}

A  bipartite  graph $G$ is $k$-fault-tolerant Hamiltonian laceable if $G-F$ remains Hamiltonian laceable for $F\subseteq V(G)\cup E(G)$ with $|F|\le k$. A bipartite  graph $G$ is $k$-edge-fault-tolerant Hamiltonian laceable if $G-F$ remains Hamiltonian laceable for $F\subseteq E(G)$ with $|F|\le k$.  Zhou et al. obtained the following result.

\begin{lem}{\rm (\cite{Qingguozhou})}\label{L2} The balanced hypercube $BH_n$ is $(2n-2)$-edge-fault-tolerant Hamiltonian laceable for $n\ge 2$.
\end{lem}

\begin{lem}{\rm (\cite{H.Lu})}\label{L3} Let $n\ge 2$. Then, $BH_n-\partial D_0$  has four components, and each component is isomorphic to $BH_{n-1}$.
\end{lem}

 The above lemma shows that one can divide $BH_n$ into four $BH_n$s by deleting $\partial D_d$ for any $d\in\{0, 1, \cdots, n-1\}$. The four components  of $BH_n-\partial D_j$  are  $BH^{j, 0}_{n-1}$, $BH^{j, 1}_{n-1}$, $BH^{j, 2}_{n-1}$, and $BH^{j, 3}_{n-1}$ for  $1\le j\le n-1$. For convenience, we use $BH^{0, 0}_{n-1}, BH^{0, 1}_{n-1}, BH^{0, 2}_{n-1}$, and $BH^{0, 3}_{n-1}$ to denote the components of $BH_n-\partial D_0$ throughout this paper.

\begin{lem}{\rm(\cite{Chengdongqin3})}\label{L4} Let $U, V$ be  two distinct bipartitions of $BH_n$, $u^1, u^2$ be two different vertices in $U$, and $v^1, v^2$ be two different vertices in $V$. Then, there are two  disjoint paths $P, Q$ such that $(1)$ $P$ joins $u^1$ and $v^1$, and $Q$ joins $u^2$ and $v^2$; $(2)$ $V(P\cup Q)=V(BH_n)$.

\end{lem}

A graph $G$ is hyper-Hamiltonian laceable if it is Hamiltonian laceable and, for an arbitrary vertex $v$ in $V_{i}$ where $i\in \{0, 1\}$, there is a Hamiltonian path in $G-v$ joining any two different vertices in $V_{1-i}$.  L\"{u} et al. obtained the following result.

\begin{lem}{\rm (\cite{Huazhonglu})}\label{L5} The balanced hypercube $BH_n$ is hyper-Hamiltonian laceable for $n\ge 1$.
\end{lem}

\section{Main result}
In this section, we will give the proof of the main result. First, we introduce some lemmas that will be used in the proof of the main result.

\begin{defin} Suppose that $G$ is a graph and $F\subseteq E(G)$. A vertex $u$ is called $i$-rescuable in $G$ if $|N_{G-F}(u)|=i$.
\end{defin}

\begin{lem}\label{P2} Let $n\ge 3, F\subseteq E(BH_n)$ with $|F|=4n-5$ and $\delta( BH_n-F)\ge 2$. Then, at least one statement holds in the following $(1)$ and $(2)$.

$(1)$ There is an integer $m$ in $\{0, 1, \cdots, n-1\}$ such that $|F\cap\partial D_m|\ge3$ and $\delta(BH_n-\partial D_m)\ge 2$.

$(2)$ There are two integers $m, m'$ in $\{0, 1, \cdots, n-1\}$ such that   $|F\cap\partial D_m|\ge2$ and $|F\cap\partial D_{m'}|\ge2$. Furthermore, there is no isolated vertex and no more than one $1$-rescuable vertex in $BH_n-\partial D_m$ (resp.,   $BH_n-\partial  D_{m'}$).
\end{lem}

\noindent{\bf Proof: } Let $F$ be a faulty set in $E(BH_n)$ with $|F|=4n-5$ and $\delta(BH_n-F)\ge 2$.
Suppose that $u$ is $k$-rescuable, and $v$ is $t$-rescuable in $BH_n$ where $k\le t\le  {\rm min}\{m\mid x$ is $m$-rescuable in $BH_n$, where $x\in V(BH_n)\setminus \{u, v\}\}$.

 If $k\ge 4$, there is an integer $m\in\{0, 1, \cdots, n-1\}$ such that $|F\cap \partial D_m|\ge 3$ because $|F|=4n-5>2n$.  Each vertex is at least $2$-rescuable in $BH_n-\partial D_m$ for $k\ge 4$. Thus,  statement $(1)$ holds. In the following, we assume that $2\le k\le 3$ and suppose that the dimensions of fault-free edges  incident with $u$ are $i_1, i_2, \cdots, i_k$ and the dimensions of fault-free edges incident with $v$ are $j_1, j_2, \cdots, j_t$ where $\{i_1, i_2, \cdots, i_k\}$ and $\{j_1, j_2, \cdots, j_t\}$ may be multisets.

\noindent {\bf Case 1: } $k=2$.

\noindent{\bf Subcase 1.1: } $t=2$.

In this subcase, $|F\cap \{(u, x)\mid x \in N_{BH_n}(u)\}|=|F\cap \{(v, x)\mid x \in N_{BH_n}(v)\}|=2n-2$.
Since $2(2n-2)-(4n-5)=1$, we obtain $(u, v)\in F$, and there is no faulty edge in $BH_n-\{u, v\}$.
Note that there is no triangle in $BH_n$. For an arbitrary vertex $x$ in $V(BH_n)\setminus \{u, v\}$, $x$ is incident with no more than one faulty edge, which means that $x$ is  at least $3$-rescuable in $BH_n-\partial D_i$ for all $i\in \{0, 1, \cdots, n-1\}$ as $2(n-1)-1\ge3$.
Note that for $4n-5>2n$, there is an $m\in \{0, 1, \cdots, n-1\}$ such that $|F\cap \partial D_m|\ge 3$.

  If $m\notin\{i_1, i_2, j_1, j_2\}$, then $u$ and $v$ are $2$-rescuable in $BH_n-\partial D_m$.  Therefore, statement $(1)$ holds.

 If $m\in \{i_1, i_2, j_1, j_2\}$,  without loss of generality, let $i_1=m$. Note that each vertex in $BH_n$ is incident with two $i$-dimension edges for $i\in \{0, 1, \cdots, n-1\}$. The integer $m$ appears no more than once in the multiset $\{i_1, i_2, j_1, j_2\}$.
 Thus, $u$ is $1$-rescuable, and $v$ is $2$-rescuable in $BH_n-\partial D_m$.
 Note that for $n\ge 3$, according to the pigeonhole principle, there is an integer $m'\in\{0, 1, \cdots, n-1\}\setminus\{m\}$ such that $m'$ appears no more than once in the multiset $\{ i_2, j_1, j_2\}$. Thus, $|F\cap\partial D_{m'}|\ge 2$. Without loss of generality, we can assume that $u$ is $1$-rescuable and $v$ is $2$-rescuable in $BH_n-\partial D_{m'}$. Therefore,  statement $(2)$ holds.

\noindent{\bf Subcase 1.2: } $t=3$.

First, we prove that  each vertex in $V(BH_n)\setminus \{u, v\}$ is at least $2$-rescuable in $BH_n-\partial D_i$ for all $i\in\{0, 1, \cdots, n-1\}$.

 If $(u, v)\notin F$, there are  $4n-5$ faulty edges  incident with $u$ or $v$. This means there is no faulty edge in $BH_n-\{u, v\}$. Thus, each vertex in $V(BH_n)\setminus \{u, v\}$ is incident with no more than two faulty edges. That means it is at least $2$-rescuable in $BH_n-\partial D_i$ for all $i\in\{0, 1, \cdots, n-1\}$ since $2(n-1)-2\ge 2$.

 If $(u, v)\in F$, there are $4n-6$ faulty edges  incident with $u$ or $v$. This  means that there is no more than one faulty edge in $BH_n-\{u, v\}$.  Since there is no triangle in $BH_n$, each vertex in $V(BH_n)\setminus \{u, v\}$ is incident with no more than  two faulty edges. That implies that it is at least $2$-rescuable in $BH_n-\partial D_i$ for all $i\in\{0, 1, \cdots, n-1\}$, as $2(n-1)-2\ge 2$.

 Next, we complete the proof of subcase $1.2$.

If there is an $m\in \{0, 1, \cdots, n-1\}$ such that $m\notin \{i_1, i_2, j_1, j_2, j_3\}$, then $|F\cap \partial D_m|\ge 3$ and $u$ and $v$ are at least $2$-rescuable in $BH_n-\partial D_m$. Therefore, statement $(1)$ holds.

Now, we assume that $\{0, 1, \cdots, n-1\}\subseteq\{i_1, i_2, j_1, j_2, j_3\}$.
If $n\ge 4$ or $n=3$ and $i_1=i_2$, there are two distinct integers $m, m'$ in $\{0, 1, \cdots, n-1\}\setminus\{i_1, i_2\}$. Note that $u$ is incident with two faulty $m$ (resp.,  $m'$)-dimension edges and $u$ is $2$-rescuable, $v$ is at least $1$-rescuable in $BH_n-\partial D_m$ (resp.,  $BH_n-\partial D_{m'}$). Therefore,  statement $(2)$ holds.

If $n=3$ and $i_1\not=i_2$, without loss of generality, let $i_1=0, i_2=1$. Then $2\in \{j_1, j_2, j_3\}$. Thus, $u$ is $2$-rescuable, and $v$ is at least $1$-rescuable in $BH_n-\partial D_2$.  Let $F'$ be the set of faulty edges  incident with $u$ or $v$.
Then $|F'|\ge 4n-6=6$ and $|F\cap \partial D_2|\ge |F'\cap \partial D_2|\ge 2$.
Moreover, $2\le |F'\cap \partial D_2|\le 3$.
Further, we obtain  $|F'\cap \partial D_0|\ge 2$ or $|F'\cap \partial D_1|\ge 2$, as $|F'|\ge4n-6=6$. Without loss of generality, we can assume that $|F'\cap \partial D_0|\ge 2$. Since $u$ is incident with exactly one faulty $0$-dimension edge,  $v$ is incident with at least one faulty $0$-dimension edge. Thus, $0$ appears no more than  once in the multiset $\{j_1, j_2, j_3\}$.
 Note that $u$ is $1$-rescuable, $v$ is at least $2$-rescuable in $BH_n-\partial D_0$ and $|F\cap \partial D_0|\ge|F'\cap\partial D_0|\ge 2$.  Let $m=0, m'=2$. Therefore, statement $(2)$ holds.

\noindent{\bf Subcase 1.3: } $t\ge 4$.

Since $t\ge 4$, each vertex in $BH_n-\{u\}$ is at least $2$-rescuable in $BH_n-\partial D_i$ for all $i\in\{0, 1, \cdots, n-1\}$.
Suppose that there are two distinct integers $m, m'\in \{0, 1, \cdots, n-1\}\setminus\{i_1, i_2\}$. Note that $u$ is incident with two faulty $m$ (resp.,  $m'$)-dimension faulty edges, and $u$ is $2$-rescuable in $BH_n-\partial D_m$ (resp.,  $BH_n-\partial D_{m'}$). Then, statement $(2)$ holds.

 If there aren't two distinct integers $m, m'\in \{0, 1, \cdots, n-1\}\setminus\{i_1, i_2\}$, then $n=3, i_1\not=i_2$ and $|F|=4n-5=7$. Without loss of generality, we can assume that $i_1=0, i_2=1$. Then, $|F\cap \partial D_2|\ge 2$.
If $|F\cap \partial D_2|\ge3$, then $u$ is $2$-rescuable in $BH_n-\partial D_2$. Therefore, statement $(1)$ holds.  If $|F\cap\partial D_2|=2$, then $u$ is $2$-rescuable in $BH_n-\partial D_2$.  According to the pigeonhole principle, $|F\cap \partial D_0|\ge 3$ or $|F\cap \partial D_1|\ge3$.
Without loss of generality, we can assume that $|F\cap \partial D_0|\ge 3$. Thus, $u$ is $1$-rescuable in $BH_n-\partial D_0$.  Let $m=0, m'=2$. Therefore,  statement $(2)$ holds.

\noindent{\bf Case 2: } $k=3$.

\noindent{\bf Subcase 2.1: } $t=3$.

(a) There is a $3$-rescuable vertex $w$ in $BH_n$ where $w \in V(BH_n)\setminus \{u, v\}$.

Suppose that the dimensions of fault-free edges  incident with $w$ are $s_1, s_2, s_3$.
Note that $3(2n-3)-(4n-5)=2n-4\ge2$ and  there is no triangle in $BH_n$. We obtain $n=3$ and $|E(BH_n[\{u, v, w\}])|=|F\cap E(BH_n[\{u, v, w\}])|=2$.
Without loss of generality, let $(u, v), (v, w)\in F$. Since $(4\times 3-5)-(6\times3-11)=0$, $F\subseteq \{e\mid e {\rm~is~incident~to~a~vertex~in~}\{u, v, w\}\}$. Hence, each vertex in $V(BH_n)\setminus\{u, v, w\}$ is incident with no more than  two faulty edges.  This means that each vertex in $V(BH_n)\setminus\{u, v, w\}$ is at least $2$-rescuable in $BH_n-\partial D_i$ for all $i\in \{0, 1, \cdots, n-1 \}$.

Note that $4n-5=7> 2n$. There is an integer $m\in\{0, 1, \cdots, n-1\}$ such that $|F\cap \partial D_m|\ge3$. Without loss of generality,  we can assume that $m=0$. Since each vertex is incident with exactly two $0$-dimension edges, $0$ appears no more than three times in the multiset $\{i_1, i_2, i_3, j_1, j_2, j_3, s_1, s_2, s_3\}$.

If $0$ does not appear twice in the multiset $\{i_1, i_2, i_3\}$, $\{j_1, j_2, j_3\}$ and  $\{s_1, s_2, s_3\}$, then $u, v, w$ are at least $2$-rescuable in $BH_n-\partial D_0$. Thus, statement $(1)$ holds.

Otherwise, without loss of generality, assume that $i_1=i_2=0, i_3=1$. Then, $u$ is $1$-rescuable, and $v$ and $w$ are at least $2$-rescuable in $BH_n-\partial D_0$. Moreover,  $u$ is $2$-rescuable in both $BH_n-\partial D_1$ and $BH_n-\partial D_2$.

  If $2$ appears four times in the multiset $\{j_1, j_2, j_3, s_1, s_2, s_3\}$,  then $|F\cap \partial D_1|\ge 2$.  Note that $v, w$ are at least $2$-rescuable in $BH_n-\partial D_1$. Then, let $m=0, m'=1$. Statement $(2)$ holds.

 If $2$ appears no more than  three times in the multiset $\{j_1, j_2, j_3, s_1, s_2, s_3\}$, then  one of the vertices in $\{v, w\}$ is at least
 $2$-rescuable, and another is at least $1$-rescuable in $BH_n-\partial D_2$. Let $m=0, m'=2$. Statement $(2)$ holds.

(b) Suppose that each vertex in $V(BH_n)\setminus \{u, v\}$ is at least $4$-rescuable in $BH_n$.

Note that each vertex in  $V(BH_n)\setminus \{u, v\}$ is at least $2$-rescuable in $BH_n-\partial D_i$ for all $i\in \{0, 1, \cdots, n-1\}$. Let $F'$ be the set of faulty edges  incident with $u$ or $v$. Then, $|F'|\ge 4n-7$.

If there is an integer $m\in \{0, 1, \cdots, n-1\}$ such that $|F'\cap \partial D_m|\ge 3$, then $m$ appears no more than  once in the multiset $\{i_1, i_2, i_3, j_1, j_2, j_3\}$. Thus, $u$ (resp.,  $v$) is at least $2$-rescuable in $BH_n-\partial D_m$.  Statement $(1)$ holds.

  Otherwise, $n=3$ and $|F'\cap\partial D_i|\le 2$ for all $i=0, 1, 2$. Then,  there are two distinct integers $m, m'\in \{0, 1, 2\}$ such that $(2)$ holds. In fact, there are two integers $m, m'\in\{0, 1, 2\}$ such that $|F'\cap \partial D_m|=2$ and $|F'\cap \partial D_{m'}|=2$, as $|F'|\ge 4n-7=5$. Thus, $m$ (resp.,  $m'$) appears no more than  twice in the multiset $\{i_1, i_2, i_3, j_1, j_2, j_3\}$. Hence, one vertex in $\{u, v\}$ is at least $2$-rescuable and another is at least $1$-rescuable in $BH_n-\partial D_m$ (resp.,  $BH_n-\partial D_{m'}$). Statement $(2)$ holds.

\noindent{\bf Subcase 2.2: } $ t\ge 4$.

Note that for $t\ge 4$, each vertex in $BH_n-\{u\}$ is at least $2$-rescuable in $BH_n-\partial D_i$ for all $i\in \{0,  1, \cdots, n-1\}$.
Since $4n-5> 2n$, there is an integer $m$ such that $|F\cap \partial D_m|\ge 3$. If $m$ appears no more than  once in the multiset $\{i_1, i_2, i_3\}$, then $u$ is at least $2$-rescuable in $BH_n-\partial D_m$. Statement $(1)$ holds.
If $m$ appears twice in the multiset $\{i_1, i_2, i_3\}$, without loss of generality, let $i_1=i_2=m$. Then, $u$ is $1$-rescuable in $BH_n-\partial D_m$ and let $m'$ be an integer in $\{0, 1, \cdots, n-1\}\setminus\{m, i_3\}$. Then $u$ is incident with two faulty $m'$-dimension edges and $u$ is $3$-rescuable in $BH_n-\partial D_{m'}$. Hence,  statement $(2)$ holds.

Owing to the above discussion, the lemma holds. $\qed$

\begin{lem} \label{P5}
Let $n\ge 3, F\subseteq E(BH_n)$ with $|F|\le 4n-5$. Then, for any vertex $u$ in $BH^{j, i}_{n-1}$, there is a fault-free $j$-dimension edge $(v, w)$ such that $(u, v)\in E(BH^{j, i}_{n-1})$. Moreover, if $u$ is $1$-rescuable in $BH_n-\partial D_j$, then $(u, v)$ can be a faulty edge.
\end{lem}

\noindent{\bf Proof: } Each vertex in $V(BH_n)$ is incident with two $i$-dimension edges for all $i\in \{0, 1, \cdots, n-1\}$. Thus, for any vertex $u\in BH^{j, i}_{n-1}$, its degree is $2n-2$ in $BH^{j, i}_{n-1}$. There is at least one fault-free $j$-dimension edge $(v, w)$ such that $(u, v)\in E(BH^{j, i}_{n-1})$ for $2(2n-2)-(4n-5)\ge 1$. If $u$ is incident with only one fault-free edge in $BH^{j, i}_{n-1}$, then $u$ is incident with $2n-3$ faulty edges in $BH^{j, i}_{n-1}$. Thus, $|F\cap \partial D_j|\le 2n-2$.  Since $2(2n-3)-(2n-2)=2n-4\ge 2$, there is at least one fault-free edge $(w, v)$ such that $(u, v)\in F$.

 $\qed$

\begin{lem}\label{L6} Let $F\subseteq E(BH_2)$ with $|F|\le 3$ and $\delta (BH_2-F)\ge 2$. Then, each edge in $BH_2-F$ lies on a fault-free Hamiltonian cycle.
\end{lem}

\noindent{\bf Proof: } The proof is rather long, so we  provide it in Appendix A. $\qed$

\begin{thm}\label{T1} Let $F\subseteq E(BH_n)$ with $|F|\le 4n-5$ and $\delta (BH_n-F)\ge 2$. Then, each edge in $BH_n-F$ lies on a fault-free Hamiltonian cycle.
\end{thm}

{\bf Proof: } We prove this theorem by induction on $n$. By Lemma \ref{L6}, the theorem holds for $n=2$. Assume that this is true for $2\le k\le n-1$. By Lemma \ref{P2}, one may partite $BH_n$ along dimension $j$, $0\le j\le n-1$, into four $BH_{n-1}s$, denoted by $BH^{j, 0}_{n-1}, BH^{j, 1}_{n-1}, BH^{j, 2}_{n-1}$ and  $BH^{j, 3}_{n-1}$, such that $|F\cap \partial D_j|\ge 2$ and there is no isolated vertex and no more than  one $1$-rescuable vertex in $BH_n-\partial D_j$. Let $F^i=F\cap BH^{j, i}_{n-1}$ and $e=(u, v)$ be an arbitrary  edge in $BH_n-F$. We must show that there is a Hamiltonian cycle in $BH_n-F$ that contains $e$.

\noindent {\bf Case 1: } $e\in E(BH_n)-\partial D_j$.

Without loss of generality, we can assume that $e\in BH^{j, 0}_{n-1}$.

\noindent {\bf Subcase 1.1: } $|F^i|\le 4n-9$ for all $i=0, 1, 2, 3.$

\noindent {\bf Subcase 1.1.1:} $\delta(BH_n-F-\partial D_j)=1$.

By Lemma \ref{P2}, there is exactly one $1$-rescuable vertex in $BH_n-\partial D_j$, say, $w$.
Without loss of generality, we can assume that $w$ is a white vertex.

\noindent {\bf Subcase 1.1.1.1: } $w\in BH^{j, 0}_{n-1}$.

\begin{figure}[!htbp]
  \centering
  \includegraphics[width=0.4\textwidth]{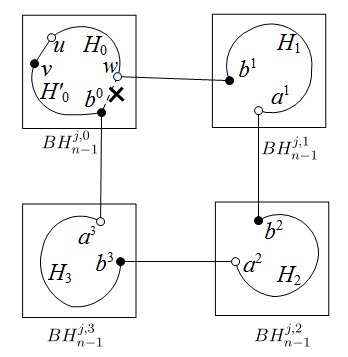}\\
  \caption{Illustration for  subcase 1.1.1.1  of theorem {\ref{T1}} }\label{1-1-1-1}
\end{figure}

 Since $\delta (BH_n-F)\ge 2$, $w$ is incident with at least one fault-free $j$ dimension edge, say, $(w, b^1)$.
 By Lemma \ref{P5}, there is a fault-free $j$-dimension edge $(b^0, a^3)$ such that $(w, b^0)\in F$.
 By induction, there is a  Hamiltonian cycle $C_0$ in $BH^{j, 0}_{n-1}-F+(w, b^0)$ that contains $(u, v)$.
 Since $w$ is incident with exactly two edges in $BH^{j, 0}_{n-1}-F+(w, b^0)$,  $(w, b^0)\in E(C_0)$.
 We represent $C_0$ as $\langle u, H_0, w, b^0, H'_0, u\rangle$.
 Note that $w$ is incident with $2n-3$ faulty edges in $BH^{j, 0}_{n-1}$, so we have $|F^0|\ge 2n-3$. Hence $|F^1|+|F^2|+|F^3|=|F|-|F^0|-|F\cap \partial D_j|\le 4n-5-(2n-3)-2=2n-4$.  Let $(a^1, b^2), (a^2, b^3)$ be fault-free $j$-dimension edges. By Lemma \ref{L2}, there is a fault-free Hamiltonian path $H_i$ in $BH^{j, i}_{n-1}$ that  joins $a^i$ and $b^i$ for $i=1, 2, 3$. Then,  the cycle $C=\langle u, H_0, w, b^1, H_1, a^1, b^2, H_2, a^2, b^3, H_3, a^3, b^0, H'_0, u\rangle$ (see figure \ref{1-1-1-1}) is the desired cycle.

\noindent {\bf Subcase 1.1.1.2: } $w\in BH^{j, 1}_{n-1}$ (or $BH^{j, 3}_{n-1}$).

 Since $4n-5-(2n-3)-2=2n-4$,  $|F^i|\le 2n-4$ for all $i=0, 2, 3$.

 Let $X=\{x^0\mid (x^0, b^1)\in \partial D_j\setminus F$, where $x^0\in BH^{j, 0}_{n-1}$ and $(w, b^1)\in  F$\}.  Since $2(2n-3)-[(4n-5)-(2n-3)]=2n-4\ge 2$, we  can obtain  $|X|\ge 1$. Moreover, if $|X|=1$, there are $2n-2$ faulty edges between $BH^{j, 0}_{n-1}$ and $BH^{j, 1}_{n-1}$.

\noindent{\bf Subcase 1.1.1.2.1: } $X=\{u\}$.

\begin{figure}[!htbp]
  \centering
  \includegraphics[width=0.8\textwidth]{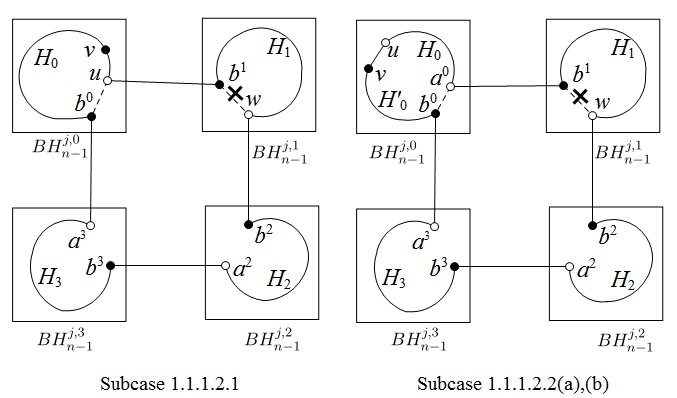}\\
  \caption{Illustration for Subcase 1.1.1.2 of theorem {\ref{T1}} }\label{1-1-1-2}
\end{figure}

 By induction, there is a fault-free Hamiltonian cycle $C_0$ in $BH^{j, 0}_{n-1}$ that contains $(u, v)$ and a Hamiltonian cycle $C_1$ in $BH^{j, 1}_{n-1}-F+(w, b^1)$ that contains $(w, b^1)$.
 Suppose that $N_{C_0}(u)=\{v, b^0\}$. Let $H_0=C_0-(u, b^0), H_1=C_1-(w, b^1)$.
Note that $|X|=1$, so there are $2n-2$ faulty $j$-dimension edges between $BH^{j, 0}_{n-1}$ and $BH^{j, 1}_{n-1}$.
Thus, $F\cap E(BH^{j, 0}_{n-1}, BH^{j, 3}_{n-1})=\emptyset$. Let $(b^0, a^3), (a^2, b^3)$ be  fault-free $j$-dimension edges.
Note that $|F^i|\le 2n-4$ for $i=2, 3$.
By Lemma \ref{L2}, there is a Hamiltonian path $H_i$ of $BH^{j, i}_{n-1}$ that joins $a^i$ and $b^i$ for $i=2, 3$. Hence, the cycle $C=\langle u, b^1, H_1, w, b^2, H_2, a^2, b^3, H_3, a^3, b^0, H_0, u\rangle $ (see figure \ref{1-1-1-2}) is the desired cycle.

\noindent{\bf Subcase 1.1.1.2.2: }  $X\not= \{u\}$.

(a) $|F\cap \partial D_j|\le 3$.

Let $a^0\in X\setminus\{u\}$ and $(w, b^2)$ be a fault-free $j$-dimension edge.
By induction, there is a fault-free Hamiltonian cycle $C_0$ in $BH^{j, 0}_{n-1}$ that contains $(u, v)$.
Suppose that $N_{C_0}(a^0)=\{b^0, d^0\}$.
Note $|F\cap \partial D_j|\le 3$; we can see that $b^0$ or $d^0$ is incident with one fault-free $j$-dimension edge.
Without loss of generality, we can assume that $b^0$ is incident with one fault-free $j$-dimension edge, say, $(b^0, a^3)$.
We can represent $C_0$ as $\langle u, H_0, a^0, b^0, H'_0, u\rangle$.
Let $(a^2, b^3)$ be a fault-free $j$-dimension edge. Note that $|F^i|\le (4n-5)-(2n-3)-2= 2n-4$ for $i=2, 3$.
By Lemma \ref{L2}, there is a fault-free Hamiltonian path $H_i$ in $BH^{j, i}_{n-1}$ that joins $a^i$ and $b^i$ for $i=2, 3$.
By induction, there is a fault-free Hamiltonian cycle $C_1$ in $BH^{j, 1}_{n-1}$ that contains $(w, b^1)$ as $|F^1|\le 4n-9$ and $\delta (BH^{j, 1}_{n-1}-F+(w, b^1))\ge 2$.
Let $H_1=C_1-(w, b^1)$.
Then, the cycle $C=\langle u, H_0, a^0, b^1, H_1, w, b^2, H_2, a^2, b^3, H_3, a^3, b^0, H'_0, u\rangle$ (see figure \ref{1-1-1-2}) is the desired cycle.

(b)  $|F\cap \partial D_j|\ge 4$.

 Since $[2(n-1)-2]+4+(2n-3)=4n-3> 4n-5$, each vertex  is at least $3$-rescuable in $BH^{j, 0}_{n-1}$. Let $(u, \alpha)$ be a fault-free edge in $BH^{j, 0}_{n-1}$ where $\alpha\not=v$.  Let $(F^0)'=F^0\cup T$ where $T=\{(u, x)\mid x\in N_{BH^{j, 0}(u)_{n-1}}-\{v, \alpha\}\}$. Then, $|(F^0)'|\le 2n-4+[4n-5-(2n-3)-4]=4n-10$, and $\delta(BH^{j, 0}_{n-1}-F'_0)\ge 2$.

Let $a^0\in X\setminus\{u\}$. By Lemma \ref{P5}, there is  a fault-free $j$-dimension edge $(b^0, a^3)$ such that $(a^0, b^0)\in E(BH^{j, 0}_{n-1})$.
 By induction, there is a Hamiltonian cycle $C_0$ in $(BH^{j, 0}_{n-1}-(F^0)')\cup (a^0, b^0)$ that contains $(a^0, b^0)$.
 Since $u$ is incident with exactly two fault-free edges in $BH^{j, 0}_{n-1}-(F^0)'$, $e=(u, v)\in E(C_0)$.
 We represent $C_0$ as $\langle u, H_0, a^0, b^0, H'_0, u\rangle$.  Also by induction, there is a Hamiltonian cycle in $BH^{j, 1}_{n-1}-F+(w, b^1)$ that contains $(w, b^1)$.
 Let $H_1=C_1-(w, b^1)$ and  $(a^2, b^3)$ be a fault-free $j$-dimension edge.
 Note that $|F^i|\le 2n-4$ for $i=2, 3$.
 By Lemma \ref{L2}, there is a Hamiltonian path $H_i$ in $BH^{j, i}_{n-1}$ that joins $a^i$ and $b^i$ for $i=2, 3$. Then, the cycle
 $C=\langle u, H_0, a^0, b^1, H_1, w, b^2, H_2, a^2, b^3, H_3, a^3, b^0, H'_0, u\rangle$ (see figure \ref{1-1-1-2}) is the desired cycle.

 \noindent{\bf Subcase 1.1.1.3: } $w\in BH^{j, 2}_{n-1}$.

 Since $\delta(BH_n-F)\ge 2$, $w$ is incident with a  fault-free $j$-dimension edge $(w, b^3)$. By Lemma \ref{P5}, there is a fault-free $j$-dimension edge $(a^1, b^2)$ such that $(w, b^2)\in F$.
 By induction, there is a fault-free Hamiltonian cycle $C_0$ in $BH^{j, 0}_{n-1}$ that contains $(u, v)$. We represent it as $\langle c^1, c^2, \cdots, c^{2^{2n-2}}, c^1\rangle$ with $c^1=u, c^{2^{2n-2}}=v$. Let $M=\{(c^1, c^2), \cdots,  (c^{2i-1}, c^{2i}), \cdots, (c^{2^{2n-2}-1}, c^{2^{2n-2}})\}$.
 Therefore, $M$ is a set of $2^{2n-3}$ mutually disjoint edges. There is an edge $(c^{2i-1}, c^{2i})$ in $M$ such that $c^{2i-1}$ (resp.,  $c^{2i}$) is incident with a fault-free $j$-dimension edge, say, $(c^{2i-1}, b^1)$ (resp.,  $(c^{2i}, a^3)$) because $2\cdot2^{2n-3}> 4n-5$. We can represent  $C_0$ as $\langle u, H_0, c^{2i-1}, c^{2i}, H'_0, u\rangle$.

 \begin{figure}[!htbp]
  \centering
  \includegraphics[width=0.4\textwidth]{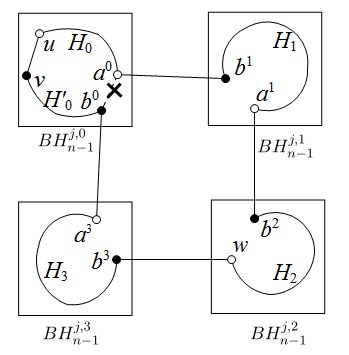}\\
  \caption{Illustration for subcase 1.1.1.3 of theorem {\ref{T1}} }\label{1-1-1-3}
\end{figure}

Note that $|F^i|\le 2n-4$ for $i=1, 3$.
 By Lemma \ref{L2}, there is a fault-free Hamiltonian path $H_i$ in $BH^{j, i}_{n-1}$ that joins $a^i$ and $b^i$ for $i=1, 3$.
 By induction, there is a Hamiltonian cycle $C_2$ in $BH^{j, 2}_{n-1}-F+(w, b^2)$ that contains $(w, b^2)$. Let $H_2=C_2-(w, b^2)$. Then, the cycle $C=\langle u, H_0, c^{2i-1}, b^1, H_1, a^1, b^2, H_2, w, b^3, H_3, a^3, c^{2i}, H'_0, u\rangle$ (see figure \ref{1-1-1-3}) is the desired cycle.

 \noindent{\bf Case 1.1.2: }  $\delta(BH_n-F-\partial D_j)\ge 2$.

\begin{figure}[!htbp]
  \centering
  \includegraphics[width=0.4\textwidth]{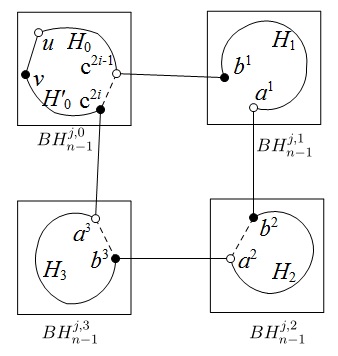}\\
  \caption{Illustration for subcase 1.1.2 of theorem {\ref{T1}} }\label{1-1-2}
\end{figure}

Since $3(2n-4)-(4n-7)=2n-5>0$, there is an integer $i\in \{1, 2, 3\}$ such that $|F^i|\le 2n-4$.
Without loss of generality, we can assume that $|F^1|\le 2n-4$.
By induction, there is a fault-free Hamiltonian cycle $C_0$ in $BH^{j, 0}_{n-1}$ that contains $(u, v)$.
Similar to the analysis of subcase 1.1.1.3, we can represent $C_0$ as $\langle u, H_0, c^{2i-1}, c^{2i}, H'_0, u\rangle$ where $(c^{2i-1}, b^1), (c^{2i}, a^{3})$ are fault-free $j$-dimension edges.
By Lemma \ref{P5}, there is  a fault-free $j$-dimension edge $(b^3, a^2)$ such that $(b^3, a^3)\in BH^{j, 3}_{n-1}$. In addition, by Lemma \ref{P5}, there is
 a fault-free $j$-dimension edge $(b^2, a^1)$ such that $(a^2, b^2)\in E(BH^{j, 2}_{n-1})$.
By induction, there is a Hamiltonian cycle $C_i$ in $(BH^{j, i}_{n-1}-F)\cup (a^i, b^i)$ that contains $(a^i, b^i)$ for $i=2, 3$. Let $H_i=C_i-(a^i, b^i)$ for i=2, 3. By Lemma \ref{L2}, there is a fault-free Hamiltonian path $H_1$ in $BH^{j, 1}_{n-1}$ that joins $a^1$ and $b^1$.
Hence, the cycle $\langle u, H_0, c^{2i-1}, b^1, H_1, a^1, b^2, H_2, a^2, b^3, H_3, a^3, c^{2i}, H'_0, u\rangle$ (see figure \ref{1-1-2}) is the desired cycle.

\noindent{\bf Subcase 1.2: } There is an integer $i\in \{0, 1, 2, 3\}$ such that $|F^i|=4n-8$.

  If $\delta (BH^{j, i}_{n-1}-F^i)\ge 2$, there are at least two vertex-disjoint  faulty edges for $(4n-8)-[2(n-1)-2]=2n-4\ge 1$. Note that $|F\cap \partial D_j|\le 3$. Then, there is a faulty edge $(u^i, v^i)$ such that $u^i$ is incident with a fault-free $j$-dimension edge and $v^i$ is incident with a fault-free $j$-dimension edge.

  Suppose that $\delta (BH^{j, i}_{n-1}-F^i)=1$. Note that there is no more than  one  $1$-rescuable vertex in $BH^{j,i}_{n-1}$. Let $u^i$ be the $1$-rescuable vertex in $BH^{j, i}_{n-1}$.
Since $\delta(BH_n-F)\ge 2$, $u^i$ is incident with a fault-free $j$-dimension edge. By Lemma \ref{P5}, there is a fault-free $j$-dimension edge $(v^i, x)$ such that $(u^i, v^i)\in F^i$. Obviously, $\delta(BH^{j, i}_{n-1}-F^i+(u^i, v^i))\ge 2$.

Owing to the above discussion,  there is an edge $(u^i, v^i)\in F^i$ such that
$(1)$ $u^i$ is incident with a fault-free $j$-dimension edge;
$(2)$ $v^i$ is incident with a fault-free $j$-dimension edge;
$(3)$ $\delta(BH^{j, i}_{n-1}-F^i+(u^i, v^i))\ge 2$, where $|F^i|=4n-8$ for  $i\in \{0, 1, 2, 3\}$.

\noindent{\bf Subcase 1.2.1: } $|F^0|=4n-8$.

\begin{figure}[!htbp]
  \centering
  \includegraphics[width=0.38\textwidth]{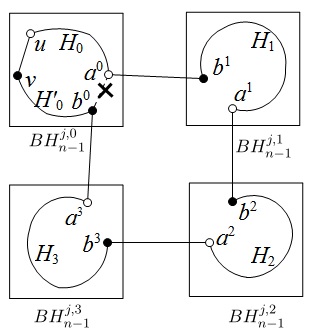}\\
  \caption{Illustration for subcase 1.2.1 of theorem {\ref{T1}} }\label{1-2-1}
\end{figure}

By induction, there is a Hamiltonian cycle $C_0$ that contains $(u, v)$ in $BH^{j, 0}_{n-1}-F^0+(u^0, v^0)$. Obviously, $|F\cap E(C_0)|\le 1$.
If $F\cap E(C_0)=1$, let $(a^0, b^0)=(u^0, v^0)$.
If $|F\cap E(C_0)|=0$, let $(a^0, b^0)$ be an edge such that $a^0$ (resp.,  $b^0$) is incident with a fault-free $j$-dimension edge.
We can represent the cycle $C_0$ as $\langle u, H_0, a^0, b^0, H'_0, u\rangle$ for $2\cdot 2^{2n-3}>4n-5$, $(a^0, b^0)$ exists.
Let $(a^0, b^1), (a^1, b^2), (a^2, b^3), (a^3, b^0)$ be fault-free $j$-dimension edges. By Lemma \ref{L2}, there is a fault-free Hamiltonian path $H_i$ that joins $(a^i, b^i)$ in $BH^{j, i}_{n-1}$ for $i=1, 2, 3$. The cycle $C=\langle u, H_0, a^0, b^1, H_1, a^1, b^2, H_2, a^2, b^3, H_3, a^3, b^0, H'_0, u\rangle$ (see figure \ref{1-2-1}) is the desired cycle.

\noindent{\bf Subcase 1.2.2: } $|F^1|=4n-8$ (or $|F^3|=4n-8$).

\begin{figure}[!htbp]
  \centering
  \includegraphics[width=0.4\textwidth]{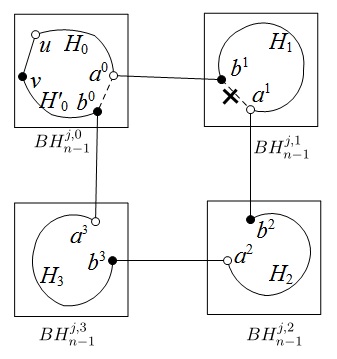}\\
  \caption{Illustration for subcase 1.2.2 of theorem {\ref{T1}} }\label{1-2-2}
\end{figure}

In this subcase, $|F^0|\le 1$. By induction, there is a  Hamiltonian cycle $C_1$ in $BH^{j, 1}_{n-1}-F^1+(u^1, v^1)$ that contains $(u^1, v^1)$, where $u^1$ (resp.,  $v^1$) is incident with a fault-free $j$-dimension edge. We can assume that $(u^1, v^1)=(a^1, b^1)$ and $(b^1, a^0)$, $(a^1, b^2)$ are fault-free $j$-dimension edges  incident with $a^1$ or $b^1$.
Let $H_1=C_1-(a^1, b^1)$. Since $|F^0|\le 1$ and $|F\cap\partial D_j|\le 3-|F^0|$, we can
  choose two fault-free edges $(a^0, b^0), (a^0, d^0)\in BH^{j, 0}_{n-1}-F^0$ such that $b^0$ is incident with a fault-free $j$-dimension edge $(b^0, a^3)$.
  Particularly, if $a^0=u$, let $d^0=v$.
  Let $(F^0)'=F^0\cup T$ where $T=\{(a^0, x)\mid x\in N_{BH^{j, 0}_{n-1}}(a^0)\setminus\{b^0, d^0\}\}$.
     Then, $|(F^0)'|\le 2n-3\le 4n-9$ and $\delta(BH^{j, 0}_{n-1}-(F^0)')\ge 2$. By induction, there is a Hamiltonian cycle $C_0$ in $BH^{j, 0}_{n-1}$ that contains $(u, v)$.
 Since $a^0$ is incident with exactly two fault-free edges in $BH^{j, 0}_{n-1}-(F^0)'$, we have $(a^0, b^0)\in E(C_0)$. We represent $C_0$ as $\langle u, H_0, a^0, b^0, H'_0, u\rangle$. Let $(a^2, b^3)$ be the fault-free $j$-dimension edge.
 By Lemma \ref{L2}, there is a fault-free Hamiltonian path $H_i$ in $BH^{j, i}_{n-1}$ that joins $a^i$ and $b^i$ for $i=2, 3$. Hence, the cycle $C=\langle u, H_0, a^0, b^1, H_1, a^1, b^2, H_2, a^2, b^3, H_3, a^3, b^0, H'_0, u\rangle$ (see figure \ref{1-2-2}) is the desired cycle.

 \noindent {\bf Subcase 1.2.3: } $|F^2|=4n-8$.

\begin{figure}[!htbp]
  \centering
    \includegraphics[width=0.38\textwidth]{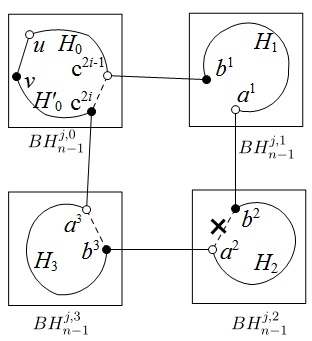}\\
  \caption{Illustration for subcase 1.2.3 of theorem {\ref{T1}} }\label{1-2-3}
\end{figure}

 By induction, there is a Hamiltonian in $BH^{j, 2}_{n-1}-F^2+(u^2, v^2)$ that contains $(u^2, v^2)$ and $u^2$ (resp.,  $v^2$) is incident with a fault-free $j$-dimension edge. We can assume that $(u^2, v^2)=(a^2, b^2)$ and $(a^2, b^3)$, $(a^1, b^2)$ are fault-free $j$-dimension edges.   Let $H_2= C_2-(a^2, b^2)$.
 By induction, there is a fault-free Hamiltonian cycle $C_0$ in $BH^{j, 0}_{n-1}$ that contains $(u, v)$. Similar to the analysis of subcase 1.1.1.3, we can represent $C_0$ as $\langle u, H_0, c^{2i-1}, c^{2i}, H'_0, u\rangle$ where $c^{2i-1}$ (resp.,  $c^{2i}$) is incident with a fault-free $j$-dimension edge $(c^{2i-1}, b^1)$ (resp.,  $(c^{2i}, b^3)$). By Lemma \ref{L2}, there is a fault-free Hamiltonian path $H_i$ of $BH^{j, i}_{n-1}$ that contains $(a^i, b^i)$ for $i=1, 3$.
 Hence, the cycle $C=\langle u, H_0, c^{2i-1}, b^1, H_1, a^1, b^2, H_2, a^2, b^3, H^3, a^3, c^{2i}$, $H'_0, u\rangle$ (see figure \ref{1-2-3}) is the desired cycle.

 \noindent{\bf Subcase 1.3: } There is an integer $i$ in $\{0, 1, 2, 3\}$ such that $|F^i|=4n-7$.

We prove this subcase according to two cases below.

\noindent{\bf Subcase 1.3.1} $\delta(BH^{j, i}_{n-1}-F^i)\ge 2$ for all $i\in\{0, 1, 2, 3\}$.

Since $4n-7-2(2n-4)=1$, there are at least three disjoint faulty edges.  There are
two faulty edges $(a^i, b^i), (c^i, d^i)$ such that $a^i, b^i, c^i, d^i$ are incident with four disjoint fault-free $j$-dimension edges owing to $|F\cap \partial D_j|=2$.

\noindent{\bf Subcase 1.3.1.1: } $|F^0|=4n-7$.

\begin{figure}[!htbp]
  \centering
  \includegraphics[width=0.95\textwidth]{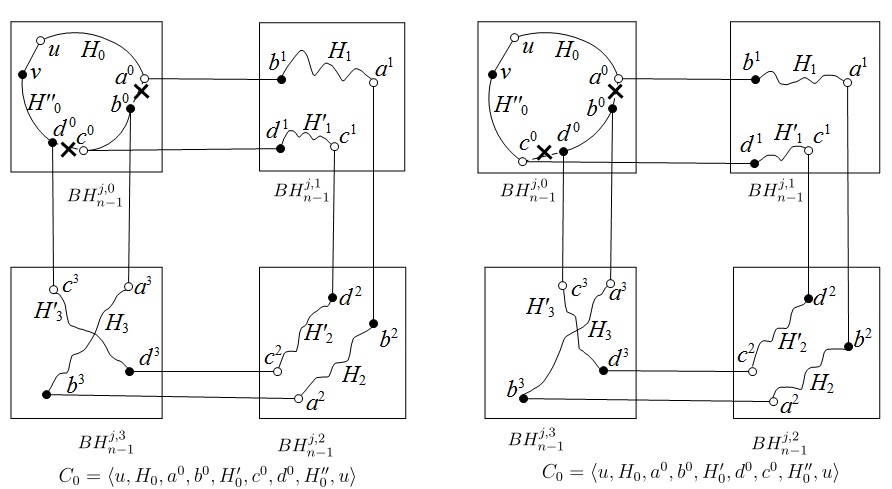}\\
  \caption{Illustration for  subcase 1.3.1.1 of theorem {\ref{T1}} }\label{1-3-1-1}
\end{figure}

By induction, there is a Hamiltonian cycle $C_0$ in $BH^{j, 0}_{n-1}-F+(a^0, b^0)+(c^0, d^0)$ that contains $(u, v)$. If $|C_0\cap \{(a^0, b^0), (c^0, d^0)\}|\le 1$, according to subcase $1.2.1$, there is a  desired cycle. If $|C_0\cap \{(a^0, b^0), (c^0, d^0)\}|=2$,  represent $C_0$ as $\langle u, H_0, a^0, b^0, H'_0, c^0, d^0, H''_0, u\rangle$ (or $\langle u, H_0, a^0, b^0, H'_0, d^0, c^0$, $H''_0, u\rangle$).  Let $(a^i, b^{i+1}), (c^i, d^{i+1})$  be fault-free $j$-dimension edges such that $a^i\not=c^i, b^i\not=d^i$ for all $i=0, 1, 2, 3$. By Lemma \ref{L4}, there are two disjoint paths $H_i, H'_i$ such that
$(1) $ $H_i$ joins $a^i$ and $b^i$, $H'_i$ joins $c^i$ and $d^i$; $(2)$ $V(H_i\cup H'_i)=V(BH^{j, i}_{n-1})$ for $i=1, 2, 3$. Hence, the cycle $C=\langle u, H_0, a^0, b^1, H_1, a^1, b^2, H_2, a^2$, $b^3, H_3, a^3, b^0, H'_0, c^0$, $d^1, H'_1, c^1, d^2, H'_2, c^2, d^3, H'_3,c^3$, $d^0, H''_0, u\rangle$(or the cycle $C=\langle u, H_0, a^0, b^1, H_1, a^1$, $b^2, H_2, a^2, b^3, H_3, a^3, b^0$, $H'_0, d^0, c^3, H'_3, d^3, c^2, H'_2, d^2, c^1$, $H'_1, d^1, c^0, H''_0, u\rangle$ ) (see figure \ref{1-3-1-1}) is the desired cycle.

\noindent{\bf Subcase 1.3.1.2: } $|F^1|=4n-7$ (or $|F^3|=4n-7$).

\begin{figure}[!htbp]
  \centering
  \includegraphics[width=0.9\textwidth]{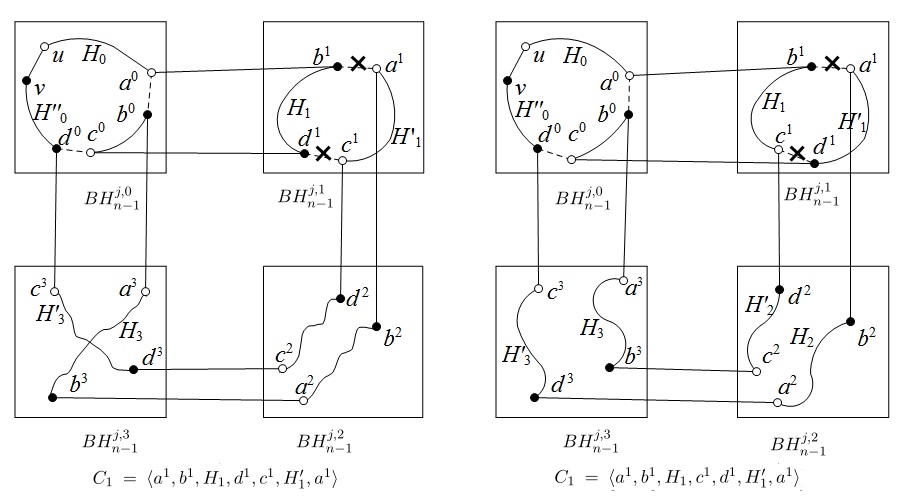}\\
  \caption{Illustration for  subcase 1.3.1.2 of theorem {\ref{T1}} }\label{1-3-1-2}
\end{figure}

By induction, there is a Hamiltonian cycle $C_1$ in $BH^{j, 1}_{n-1}-F+(a^1, b^1)+(c^1, d^1)$ that contains $(a^1, b^1)$.
If $(c^1, d^1)\notin C_1$, according to subcase $1.2.2$, there is a desired cycle. Thus,  assume that $(c^1, d^1)\in C_1$. We can assume that $C_1=\langle a^1, b^1, H_1, d^1, c^1, H'_1, a^1\rangle$ (or $C_1=\langle a^1, b^1, H_1, c^1, d^1, H'_1, a^1\rangle$) and that $(a^0, b^1), (c^0, d^1), (a^1, b^2), (c^1, d^2)$ are  fault-free $j$-dimension edges.
Note that $a^0\not=c^0$. Then, $a^0\not=u$ or $c^0\not=u$. Without loss of generality, we assume that $c^0\not=u$.
If there is a vertex $x^0$ in $BH^{j, 0}_{n-1}$  incident with two faulty $j$-dimension edges, let $(F^0)'=E(BH_n)\cap\{(a^0, x^0), (c^0, x^0)\}$.  Otherwise, let $(F^0)'=\{(c^0, \alpha)\mid \alpha\in BH^{j, 0}_{n-1}, {\rm ~where~} \alpha {\rm ~is~ incdent~ with~ a~ faulty~}$ $j\rm{-dimension~ edge}\}$. Since $|F\cap \partial D_j|= 2$, $|(F^0)'|\le 2$.
By induction, there is a Hamiltonian cycle $C_0$ in $BH^{j, 0}_{n-1}-(F^0)'$ that contains $(u, v)$. We may represent it as $\langle u, H_0, a^0, b^0, H'_0, c^0, d^0, H''_0, u\rangle$.  Let $(b^0, a^3), (d^0, c^3), (a^2, b^3), (c^2, d^3)$ be fault-free $j$-dimension edges where $a^2\not=c^2, b^3\not=d^3 , a^3\not=c^3$ (when $C_1=\langle a^1, b^1, H_1, c^1, d^1, H'_1, a^1\rangle$, let $(b^0, a^3), (d^0, c^3), (a^2, d^3), (c^2, b^3)$ be fault-free $j$-dimension edges).  By Lemma \ref{L4}, there are two disjoint paths $H_i, H'_i$ such that
$(1) $ $H_i$ joins $a^i$ and $b^i$, and $H'_i$ joins $c^i$ and $d^i$;  $(2)$ $V(H_i\cup H'_i)=V(BH^{j, i}_{n-1})$ for $i=2, 3$. Hence, The cycle $C=\langle u, H_0, a^0, b^1, H_1, d^1, c^0, H'_0, b^0, a^3, H_3, b^3, a^2$, $H_2, b^2, a^1, H'_1, c^1, d^2, H'_2, c^2,$ $d^3, H'_3, c^3, d^0, H''_0, u\rangle$ (or $C=\langle u, H_0, a^0, b^1, H_1, c^1, d^2, H'_2, c^2, b^3, H_3, a^3, b^0, H'_0$, $c^0, d^1, H'_1, a^1, b^2, H_2, a^2$, $d^3, H'_3, c^3, d^0, H''_0, u\rangle$) (see figure \ref{1-3-1-2}) is the desired cycle.

\noindent{\bf Subcase 1.3.1.3: } $|F^2|\le 4n-7$.
 \begin{figure}[!htbp]
  \centering
  \includegraphics[width=0.89\textwidth]{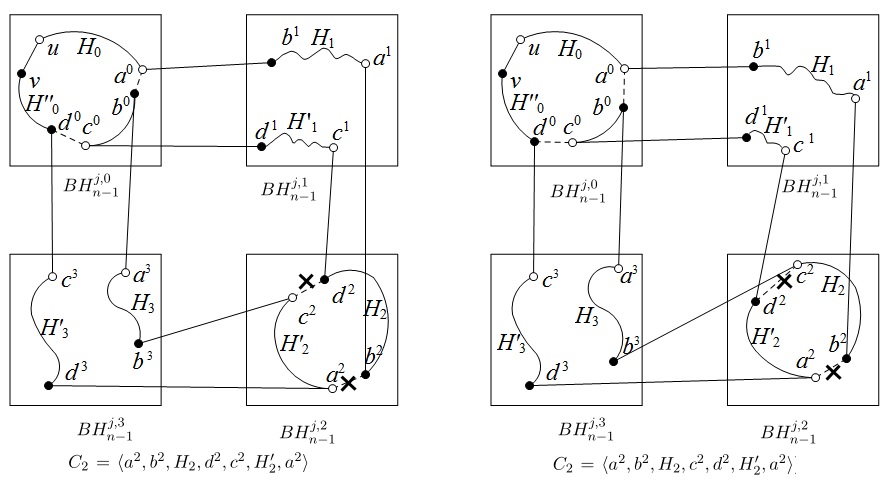}\\
  \caption{Illustration for  subcase 1.3.1.3 of theorem {\ref{T1}} }\label{1-3-1-3}
\end{figure}

By induction, there is a Hamiltonian cycle $C_2$ in $BH^{j, 2}_{n-1}-F+(a^2, b^2)+(c^2, d^2)$ that contains $(a^2, b^2)$. If $(c^2, d^2)\notin E(C_2)$, according to subcase $1.2.3$, there is a desired cycle. Now, let $(c^2, d^2)\in E(C_2)$. We can assume that $C_2=\langle a^2, b^2, H_2, d^2, c^2, H'_2, a^2\rangle$ (or $C_2=\langle a^2, b^2, H_2, c^2, d^2, H'_2, a^2\rangle$).
Let $(a^1, b^2), (c^1, d^2), (a^2, d^3), (c^2, b^3)$ be fault-free $j$-dimension edges, where $a^1\not=c^1, b^3\not=d^3$.
By induction, there is a Hamiltonian cycle $C_0$ that contains $(u, v)$. Since $|F\cap \partial D_j|= 2$, there are two edges $(a^0, b^0), (c^0, d^0)\in E(C_0)$ such that $(a^0, b^1), (c^0, d^1), (b^0, a^3), (d^0, c^3)$ are fault-free $j$-dimension edges, where $b^1\not=d^1, a^3\not=c^3$. We can represent $C_0$ as $\langle u, H_0, a^0, b^0, H'_0, c^0, d^0, H''_0, u\rangle$.
By  Lemma \ref{L4}, there are two disjoint paths $H_i, H'_i$ such that
$(1) $ $H_i$ joins $a^i$ and $b^i$, and $H'_i$ joins $c^i$ and $d^i$; $(2)$ $V(H_i\cup H'_i)=V(BH^{j, i}_{n-1})$ for $i=1,  3$. Thus, the cycle $\langle u, H_0, a^0, b^1, H_1, a^1, b^2, H_2, d^2, c^1, H'_1, d^1, c^0, H'_0,$ $ b^0, a^3, H_3, b^3, c^2, H'_2, a^2, d^3, H'_3, c^3, d^0, H''_0, u\rangle$(or $C=\langle u, H_0, a^0, b^1, H_1, a^1, b^2, H_2, c^2, b^3, H_3, a^3, b^0, H'_0, c^0,$ $d^1, H'_1, c^1, d^2, H'_2, a^2, d^3, H'_3, c^3, d^0, H''_0, u\rangle$) (see figure \ref{1-3-1-3}) is the desired cycle.

\noindent{\bf Subcase 1.3.2: } $\delta(BH^{j, i}_{n-1}-F^i)=1$ for some $i\in \{0, 1, 2, 3\}$.

By Lemma \ref{P2}, there is no more than  one $1$-rescuable vertex in $BH^{j, i}_{n-1}$, say, $w$. Suppose that $(w, \beta)$ is the fault-free edge incident with $w$.  Without loss of generality, we can assume that $w$ is a white vertex. Since $2n-3\ge 3$, $w$ is incident with at least three faulty edges in $BH^{j, i}_{n-1}$. Then, there are two faulty edges $(w, b^i), (w, d^i)$ such that $b^i, d^i$ are incident with two disjoint fault-free $j$-dimension edges owing to $|F\cap\partial D_j|=2$.

\noindent{\bf Subcase 1.3.2.1: } $|F^0|=4n-7$.

 By induction, there is a Hamiltonian cycle $C_0$ in $BH^{j, 0}_{n-1}-F+(w, b^0)+(w, d^0)$ that contains $(u, v)$.
 If $|C_0\cap \{(w, b^0), (w, d^0)\}|\le 1$, according to  subcase $1.2.1$, there is a desired cycle. Now, we can assume that $|C_0\cap \{(w, b^0), (w, d^0)\}|=2$.

 (a) Suppose that $w$ is incident with two fault-free $j$-dimension edges.
 \begin{figure}[!htbp]
  \centering
  \includegraphics[width=0.4\textwidth]{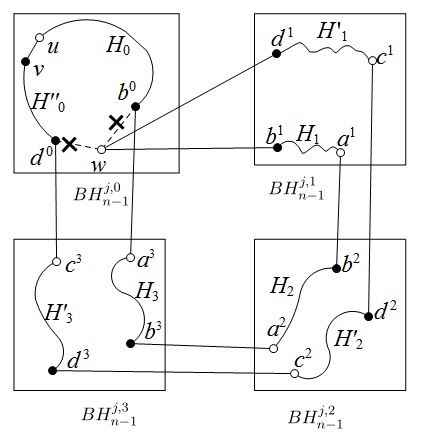}\\
  \caption{Illustration for subcase 1.3.2.1 (a) of theorem {\ref{T1}} }\label{1-3-2-1a}
\end{figure}
 In this instance, we can represent $C_0$ as $\langle u, H_0, b^0, w, d^0, H'_0, v, u\rangle$.
 Let $(w, b^1), (w, d^1)$, $(a^1, b^2), (c^1, d^2)$, $(a^2, b^3), (c^2, d^3), (a^3, b^0), (c^3, d^0)$ be fault-free $j$-dimension edges. By Lemma \ref{L4}, there are two paths $H_i$ and $H'_i$ such that $(1) H_i$ joins $a^i$ and $b^i$, and $H'_i$ joins $c^i$ and $d^i$; $(2) V(H_i\cup H'_i)=V(BH^{j, i}_{n-1})$ for $i=1, 2, 3$. Thus, the cycle $\langle u, H_0, b^0, a^3, H_3, b^3, a^2, H_2, b^2, a^1$, $H_1, b^1, w, d^1, H'_1, c^1, d^2, H'_2, c^2, d^3, H'_3, c^3, d^0, H'_0, u\rangle$ (see figure \ref{1-3-2-1a}) is the desired cycle.

(b) Suppose that $w$ is incident with exactly one fault-free $j$-dimension edge $(w, b^1)$.

\begin{figure}[!htbp]
  \centering
  \includegraphics[width=0.9\textwidth]{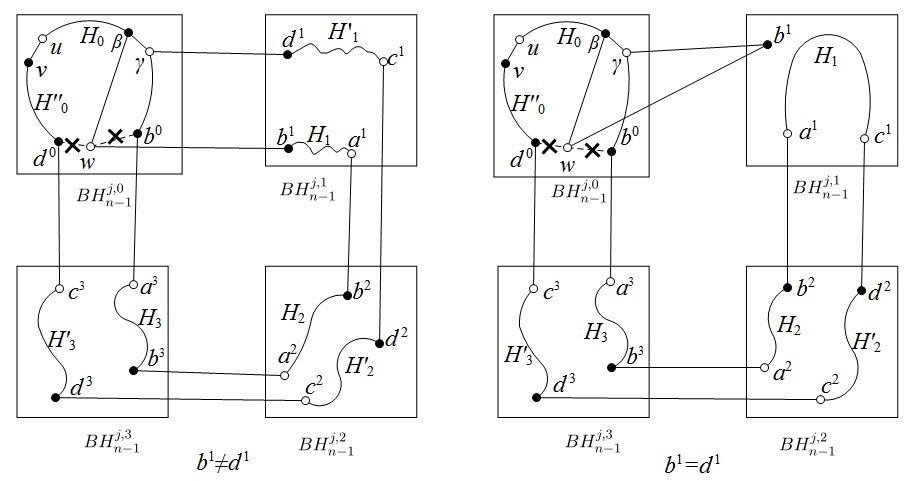}\\
  \caption{Illustration for subcase 1.3.2.1 (b) of theorem {\ref{T1}} }\label{1-3-2-1b}
\end{figure}

In this instance, we can represent $C_0$ as $\langle u, H_0, \beta, \gamma, H'_0, b^0, w, d^0, H''_0, u\rangle$.
Suppose that $\gamma$ is incident with a fault-free $j$-dimension edge $(\gamma, d^1)$ where $b^1\not=d^1$. Let $(b^0, a^3)$, $(d^0, c^3)$ and  $(a^i, b^{i+1}), (c^i, d^{i+1})$ be fault-free $j$-dimension edges for $i=1, 2$.
 By Lemma \ref{L4}, there are two paths $H_i$ and $H'_i$ such that $(1) H_i$ joins $a^i$ and $b^i$, and $H'_i$ joins $c^i$ and $d^i$; $(2) V(H_i\cup H'_i)=V(BH^{j, i}_{n-1})$ for $i=1, 2, 3$. Thus, the cycle $C=\langle u, H_0, \beta, w, b^1, H_1, a^1, b^2, H_2, a^2, b^3, H_3, a^3, b^0, H'_0, \gamma, d^1$, $H'_1, c^1, d^2, H'_2$, $c^2, d^3, H'_3, c^3, d^0, H''_0, u\rangle$ (see figure \ref{1-3-2-1b}) is the desired cycle.

 Suppose that the unique fault-free $j$-dimension edge  incident with $\gamma$ is $(\gamma, b^1)$. Let $(b^0, a^3)$, $(d^0, c^3)$ and $(a^i, b^{i+1}), (c^i, d^{i+1})$ be fault-free $j$-dimension edges for $i=1, 2$. By Lemma \ref{L4}, there are two paths $H_i$ and $H'_i$ in $BH^{j, i}_{n-1}$ such that $(1) H_i$ joins $a^i$ and $b^i$, and $H'_i$ joins $c^i$ and $d^i$; $(2) V(H_i\cup H'_i)=V(BH^{j, i}_{n-1})$ for $i=2, 3$. By Lemma \ref{L5}, there is a Hamiltonian path in $BH^{j, 1}_{n-1}-\{\alpha\}$ that joins $a^1$ and $c^1$. Then, the cycle $C=\langle u, H_0, \beta, w, b^1, \gamma, H'_0, b^0, a^3, H_3, b^3, a^2, H_2, b^2, a^1, H_1,$ $ c^1, d^2, H'_2, c^2, d^3, H'_3, c^3, d^0$, $H''_0, u\rangle$ (see figure \ref{1-3-2-1b}) is the desired cycle.

\noindent{\bf Subcase 1.3.2.2: } $|F^i|=4n-7$ for $i=1, 2, 3$.

By induction, there is a Hamiltonian cycle $C_i$ in $BH^{j, i}_{n-1}-F+(w, b^i)+(w, d^i)$ that contains $(w, \beta)$. Then, $|C^i\cap \{(w, b^i), (w, d^i)\}|\le 1$.  According to subcase $1.2.2$ or $1.2.3$, there is a desired cycle.

\noindent {\bf Case 2: } $e\in \partial D_j$.

If $|F\cap\partial D_j|=2$ or $\delta(BH_n-F-\partial D_j)=1$, by Lemma \ref{P2}, there is an integer $j'\in \{0, 1, \cdots, n-1\}\setminus\{j\}$ such that $|F\cap\partial D_{j'}|\ge2$ and there is no isolated vertex and no more than  one $1$-rescuable vertex in $BH_n-\partial D_{j'}$. Note that $e\in BH_n-\partial D_{j'}$. According to case $1$, there is a desired cycle. Thus,  we only need to show that there is a Hamiltonian cycle in $BH_n-F$ that contains $e$ when $\delta(BH_n-F-\partial D_j)\ge 2$ and $|F\cap\partial D_j|\ge 3$.
Without loss of  generality, we can assume that $u=a^0, v=b^0$.

\noindent {\bf Subcase 2.1: } $|F^i|\le 4n-9$ for all $i=0, 1, 2, 3$.

\begin{figure}[!htbp]
  \centering
  \includegraphics[width=0.35\textwidth]{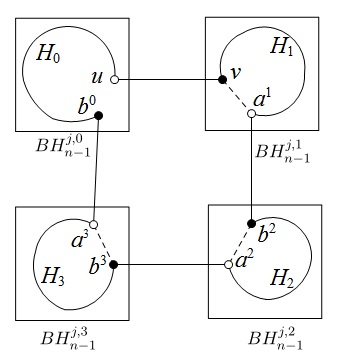}\\
  \caption{Illustration for  subcase 2.1 of theorem {\ref{T1}} }\label{2-1}
\end{figure}

Since $2(2n-3)+2> 4n-5$, we have $|F\cap BH^{j, 0}_{n-1}|\le 2n-4$ or $|F\cap BH^{j, 1}_{n-1}|\le 2n-4$. Without loss of generality, we can assume that $|F\cap BH^{j, 0}_{n-1}|\le 2n-4$. By Lemma \ref{P5}, there is a fault-free $j$-dimension edge $(a^1, b^2)$  (resp.,  ($(a^2, b^3), (a^3, b^0)$)) such that $(a^1, v)\in E(BH_n)$ (resp.,  $(a^2, b^2), (a^3, b^3)\in E(BH_n)$). By Lemma \ref{L2}, there is a Hamiltonian path $H_0$ that joins $u$ and $b^0$.
By induction, there is a Hamiltonian cycle $C_i$ in $(BH^{j, i}_{n-1}-F)\cup \{(a^i, b^i)\}$ that contains $(a^i, b^i)$ for $i=1, 2, 3$. Let $H_i=C_i-(a^i, b^i)$. The cycle $C=\langle u, v, H_1, a^1, b^2, H_2, a^2, b^3, H_3, a^3, b^0, H_0, u\rangle$ (see figure \ref{2-1}) is the desired cycle.

\noindent{\bf Subcase 2.2: } There is  an  $i\in \{0, 1, 2, 3\}$ such that $|F^i|= 4n-8$.

\noindent{\bf Subcase 2.2.1: } $|F^0|=4n-8$ (or $|F^1|=4n-8$).

 \noindent{\bf Subcase 2.2.1.1: } Suppose that $u$ is incident with at least two faulty edges in $BH^{j, 0}_{n-1}$.

\begin{figure}[!htbp]
  \centering
  \includegraphics[width=0.35\textwidth]{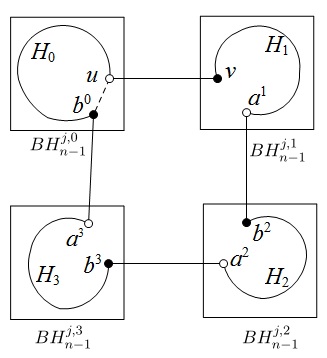}\\
  \caption{Illustration for subcase 2.2.1.1 of theorem {\ref{T1}} }\label{2-2-1-1}
\end{figure}

 Since $|F\cap \partial D_j|=3$ and each vertex is incident with exactly two $j$-dimension edges, there is a faulty edge $(u, b^0)$ such that $b^0$ is incident with a fault-free $j$-dimension edge, say, $(a^3, b^0)$. By induction, there is a Hamiltonian cycle $C_0$ in $BH^{j, 0}_{n-1}$ that contains $(u, b^0)$. Let $H_0=C_0-(u, b^0)$ and $(a^1, b^2), (a^2, b^3)$ be fault-free $j$-dimension edges. By Lemma \ref{L4}, there is a path $H_i$ in $BH^{j, i}_{n-1}$ that joins $a^i$ and $b^i$. Therefore, the cycle $\langle u, v, H_1, a^1, b^2, H_2, a^2, b^3, H_3, a^3, b^0, H_0, u\rangle$ (see figure \ref{2-2-1-1}) is the desired cycle.

\noindent{\bf Subcase 2.2.1.2:} Suppose that $u$ is incident with no more than  one faulty edge in $BH^{j, 0}_{n-1}$.

In this subcase, there are two disjoint faulty edges in $BH^{j, 0}_{n-1}-\{u\}$, say, $(c^0, d^0), (x^0, y^0)$, that are the result of $4n-8-(2n-4)-1=2n-5\ge 1$.
 Suppose that $x^0, a^0$ is a white vertex and  $M=\{z\mid (x^0, z)\in \partial D_j\setminus F$ or $(c^0, z)\in \partial D_j\setminus F$\}.

\noindent{\bf Subcase 2.2.1.2.1: }  $M\not=\{v\}$.

Without loss of generality,  we can assume that $(c^0, d^1)$ is a fault-free $j$-dimension edge where $d^1\not=v$.

(a) Suppose that $d^0$ is incident with a fault-free $j$-dimension edge $(d^0, c^3)$.

\begin{figure}[!htbp]
  \centering
  \includegraphics[width=0.8\textwidth]{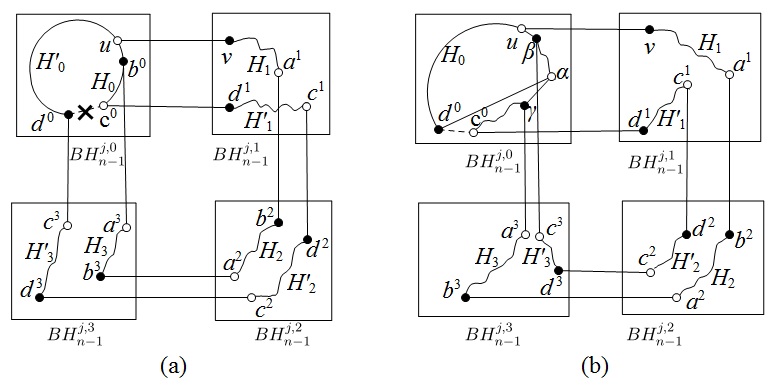}\\
  \caption{Illustration for subcase 2.2.1.2.1 (a) and (b) of theorem {\ref{T1}} }\label{2-2-1-2-1}
\end{figure}

Note that $u$ is incident with no more than  one faulty edge in $BH^{j, 0}_{n-1}$. There is a vertex $b^0\in N_{BH^{j, 0}_{n-1}}(u)$ such that $b^0$ is incident with a fault-free $j$-dimension edge $(b^0, a^3)$ where $a^3\not=c^3$ owing to $2(2n-3)-3\ge3$. By induction, there is a Hamiltonian cycle $C_0$ in $(BH^{j, 0}_{n-1}-F+(c^0, d^0))\cup \{(u, b^0)\}$ that contains $(u, b^0)$. If $(c^0, d^0)\notin E(C_0)$, similar to subcase 2.2.1.1, there is a desired cycle. Thus, we can assume that $(c^0, d^0)\in E(C_0)$.
We can represent $C_0$ as $\langle u, b^0, H_0, c^0, d^0, H'_0, u\rangle$. Let $(a^i, b^{i+1}), (c^i, d^{i+1})$ be fault-free $j$-dimension edges for $i=1, 2$ where $a^i\not=c^i, b^i\not=d^i$.
By Lemma \ref{L4}, there are two vertex-disjoint paths $H_i$ and $H'_i$ such that $(1) H_i$ joins $a^i$ and $b^i$, and $H'_i$ joins $c^i$ and $d^i$; $(2) V(H_i\cup H'_i)=V(BH^{j, i}_{n-1})$ for $i=1, 2, 3$. Hence, the cycle $C=\langle u, v, H_1, a^1, b^2, H_2, a^2, b^3, H_3, a^3, b^0,$ $H_0, c^0, d^1, H'_1, c^1, d^2, H'_2, c^2, d^3, H'_3, c^3$, $d^0, H'_0, u\rangle$ (see figure \ref{2-2-1-2-1} (a)) is the desired cycle.

(b) Suppose that $d^0$ is incident with two faulty $j$-dimension edges.

By induction, there is a Hamiltonian cycle in $BH^{j, 0}_{n-1}-F+(c^0, d^0)$ that contains $(c^0, d^0)$.
Since $\delta(BH_n-F-\partial D_j)\ge 2$, $d^0$ is incident with a fault-free edge $(d^0, \alpha)$ in $BH^{j, 0}_{n-1}$, which does  not belong to $E(C_0)$.
Without loss of generality, we can represent $C_0$ as $\langle u, H_0, d^0, c^0, H'_0,\gamma, \alpha, H''_0, \beta, u\rangle$.
Note that $d^0$ is incident with two faulty $j$-dimension edges, and $\beta, \gamma$  are incident with two disjoint fault-free $j$-dimension edges $(\beta, c^3), (\gamma, a^3)$.
Let $(a^i, b^{i+1}), (c^i, d^{i+1})$ be fault-free $j$-dimension edges where $a^i\not=c^i, b^i\not=d^i$ for $i=1, 2$. By Lemma \ref{L4}, there are two Hamiltonian paths $H_i, H'_i$ such that  $(1) H_i$ joins $a^i$ and $b^i$, and $H'_i$ joins $c^i$ and $d^i$; $(2) V(H_i\cup H'_i)=V(BH^{j, i}_{n-1})$ for $i=1, 2, 3$. Thus, the cycle $\langle u, H_0, d^0, \alpha, H''_0, \beta, c^3, H'_3, d^3, c^2, H'_2, d^2, c^1, H'_1, d^1, c^0, H'_0, \gamma, a^3, H_3, b^3, a^2, H_2, b^2$, $a^1, H_1, v, u\rangle$ (see figure \ref{2-2-1-2-1} (b)) is the desired cycle.

 \noindent{\bf Subcase 2.2.1.2.2: } $M=\{v\}$.

\begin{figure}[!htbp]
  \centering
  \includegraphics[width=0.38\textwidth]{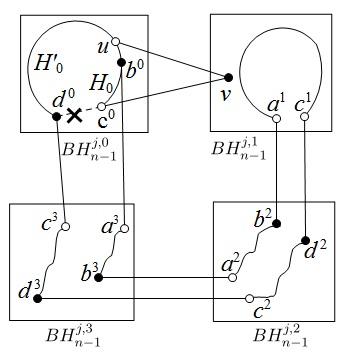}\\
  \caption{Illustration for  subcase 2.2.1.2.2 of theorem {\ref{T1}} }\label{2-2-1-2-2}
\end{figure}

In this subcase, $|E(BH^{j, 0}_{n-1}, BH^{j, 1}_{n-1})\cap F|=3$ and $|E(BH^{j, 0}_{n-1}, BH^{j, 3}_{n-1})\cap F|=0$.
By induction, there is a Hamiltonian cycle $C_0$ in $BH^{j, 0}_{n-1}-F+(c^0, d^0)$ that contains $(c^0, d^0)$. We can represent $C_0$ as $\langle u, b^0, H_0, c^0, d^0, H'_0, u\rangle$. Since $b^0$ is incident with two fault-free $j$-dimension edges, we can choose a fault-free $j$-dimension edge $(b^0, a^3)$ where $a^3\not=c^3$. Let $(a^i, b^{i+1}), (c^i, d^{i+1})$ be fault-free $j$-dimension edges for $i=1, 2$ where $a^i\not=c^i, \not=b^i\not=d^i$. By Lemma \ref{L4}, there are two vertex-disjoint paths $H_i$ and $H'_i$ such that $(1) H_i$ joins $a^i$ and $b^i$, and $H'_i$ joins $c^i$ and $d^i$; $(2) V(H_i\cup H'_i)=V(BH^{j, i}_{n-1})$ for $i= 2, 3$. By Lemma \ref{L5}, there is a Hamiltonian path $H_1$ in $BH^{j, 1}_{n-1}-\{v\}$ that joins $a^1$ and $c^1$. Hence, the cycle $C=\langle u, v, c^0, H_0, b^0, a^3, H_3, b^3, a^2, H_2, b^2, a^1, H_1, c^1, d^2, H'_2, c^2, d^3, H'_3, c^3, d^0, H'_0, u\rangle$ (see figure \ref{2-2-1-2-2}) is the desired cycle.

\noindent{\bf Subcase 2.2.2: } $|F^2|=4n-8$ (or $|F^3|=4n-8$).

\begin{figure}[!htbp]
  \centering
  \includegraphics[width=0.35\textwidth]{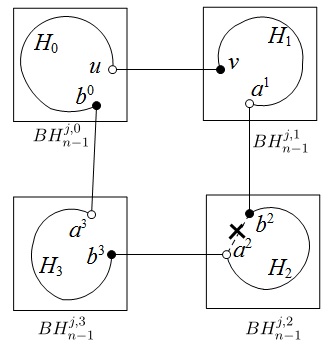}\\
  \caption{Illustration for  subcase 2.2.2 of theorem {\ref{T1}} }\label{2-2-2}
\end{figure}

There are two disjoint faulty edges in $BH^{j, 2}_{n-1}$ for $4n-8-(2n-4)=2n-4\ge 2$. Since $|F\cap\partial D_j|=3$,
there is a faulty edge $(a^2, b^2)\in BH^{j, 2}_{n-1}$ such that $a^2$ (resp.,  $b^2$) is incident with a fault-free $j$-dimension edge $(a^2, b^3)$ (resp.,  $(a^1, b^2)$). By induction there is a Hamiltonian cycle $C_2$ in $BH^{j, 2}_{n-1}-F+(a^2, b^2)$ that contains $(a^2, b^2)$. Let $H_2=C_2-(a^2, b^2)$
and  $(a^3, b^0)$ be a fault-free $j$-dimension edge. By Lemma \ref{L2}, there is  a Hamiltonian path $H_i$ in $BH^{j, i}_{n-1}$ that joins $a^i$ and $b^i$ for $i=0, 1, 3$. Thus, the cycle $C=\langle u, v, H_1, a^1, b^2, H_2, a^2, b^3, H_3, a^3, b^0, H_0, u\rangle$ (see figure \ref{2-2-2}) is the desired cycle.

\section*{Appendix A. Proof of Lemma \ref{L6}}

{\bf Lemma \ref{L6} }  Let $F\subseteq E(BH_2)$ with $|F|\le 3$ and $\delta (BH_2-F)\ge 2$. Then, each edge in $BH_2-F$ lies on a fault-free Hamiltonian cycle.

Suppose that $|F|=3$. By Lemma \ref{L1}, without loss of generality, we can assume that $|F\cap \partial D_1|\ge |F\cap \partial D_0|$.  Let $e=(u, v)$. In the following, we can assume that $a^i, c^i$ are  white vertices and $b^i, d^i$ are black vertices  in $BH^{1, i}_{1}$ for $i=0, 1, 2, 3$.

\noindent{\bf Case 1: } $|F\cap \partial D_1|=3$, $|F\cap \partial D_0|=0$.

\noindent{\bf Subcase 1.1: } $e\in \partial D_0$.

\begin{figure}[!htbp]
  \centering
  \includegraphics[width=0.8\textwidth]{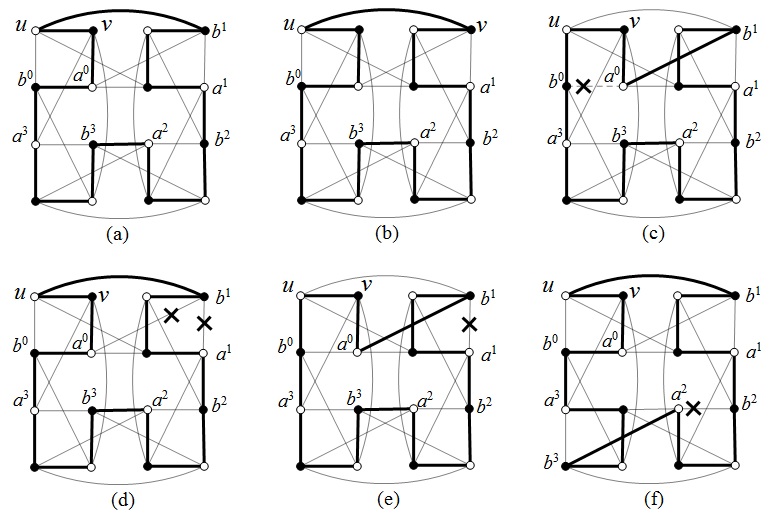}\\
   \includegraphics[width=0.8\textwidth]{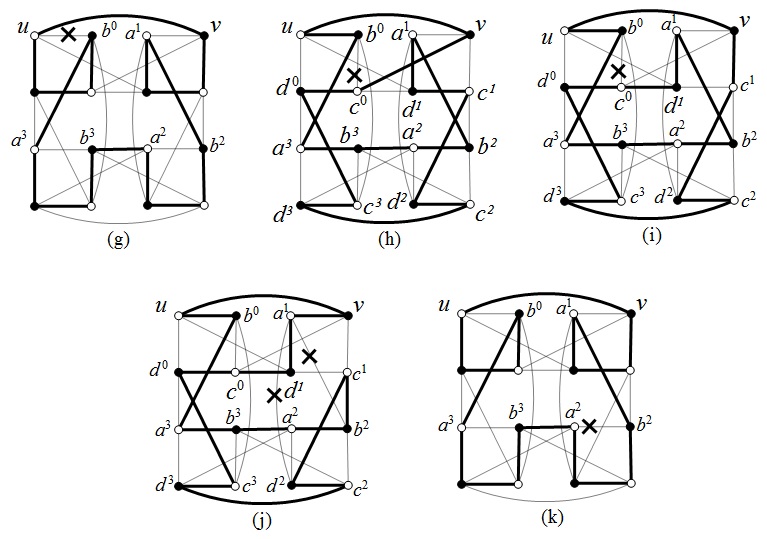}\\
  \caption{Illustration for  Lemma \ref{L6}}\label{appendix}
\end{figure}

Without loss of generality, we can assume that $e\in BH^{1, 0}_{1}$.  Since $|F\cap \partial D_1|=3$, there is an edge $(x, y)$ in $\{(v, a^0), (u, b^0)\}$ such that $x$ (resp.,  $y$) is incident with a  fault-free $1$-dimension edge. Without loss of generality, we can assume that $(x, y)=(u, b^0)$.
Let $H_0=BH^{1, 0}_1-(u, b^0)$ and $(u, b^1), (b^0, a^3), (a^1, b^2), (a^2, b^3)$ be fault-free $1$-dimension edges. Note that $|F\cap \partial D_0|=0$.
There is a Hamiltonian path $H_i$ in $BH^{1, i}_1$ that joins $a^i$ and $b^i$ for $i=1, 2, 3$.  Thus, the cycle $C=\langle u, b^1, H_1, a^1, b^2, H_2, a^2, b^3, H_3, a^3, b_0,  H_0, u\rangle$ (see figure \ref{appendix})(a) is the desired cycle.

\noindent{\bf Subcase 1.2: } $e\in \partial D_1$.

Without loss of generality, we can assume that $(u, v)$ is an edge between $BH^{1, 0}_1$ and $BH^{1, 1}_1$ and $u=a^0, v=b^1$. Let $(a^1, b^2), (a^2, b^3), (a^3, b^0)$ be fault-free $1$-dimension edges. Since $|F\cap \partial D_0|=0$, there is a Hamiltonian path $H_i$ of $BH^{1, i}_{1}$ that joins $a^i$ and $b^i$ for $i=0, 1, 2, 3$. Thus, the cycle $\langle u, v, H_1, a^1, b^2, H_2, a^2, b^3, H_3, a^3, b^0, H_0, u\rangle$ (see figure \ref{appendix} (b)) is the desired cycle.

\noindent{\bf Case 2: } $|F\cap \partial D_1|=2$, $|F\cap \partial D_0|=1$.

\noindent{\bf Subcase 2.1: } $e\in \partial D_0$.

 Without loss of generality, we can assume that $e\in BH^{1, 0}_{1}$.

 \noindent{\bf Subcase 2.1.1: } $|F^0|=1$.

 Suppose that $(a^0, b^0)$ is a faulty edge in $BH^{1, 0}_1$.  Let $H_0=BH^{1, 0}_1-(a^0, b^0)$. Then, $H_0$ is a fault-free Hamiltonian path in $BH^{1, 0}_1$ that contains $(u, v)$. Since $\delta (BH_2-F)\ge 2$, $a^0$ (resp.,  $b^0$) is incident with a fault-free $1$-dimension edge $(a^0, b^1)$ (resp.,  $(b^0, a^3)$). Let $(a^1, b^2), (a^2, b^3)$ be fault-free $1$-dimension edges. Since $|F^1|=|F^2|=|F^3|=0$, there is a Hamiltonian path $H_i$ in $BH^{1, i}_1$ that joins $a^i$ and $b^i$ for $i=1, 2, 3$. Thus, the cycle $\langle a^0, b^1, H_1, a^1, b^2, H_2, a^2, b^3, H_3, a^3, b^0, H_0, a^0\rangle$ (see figure \ref{appendix} (c))is the desired cycle.

 \noindent{\bf Subcase 2.1.2: } $|F^1|=1$(or $|F^3|=1$).

 Let $V(BH^{1, 0}_1)=\{u, v, a^0, b^0\}$ and $(a^1, b^1)$ be the faulty edge in $BH^{1, 1}_1$. Suppose that $H_1=BH^{1, 1}_1-(a^1, b^1)$. Then, $H_1$ is a Hamiltonian path in $BH^{1, 1}_1$ that joins $a^1$ and $b^1$.  Since $\delta(BH_2-F)\ge 2$, $a^1$ (resp.,  $b^1$) is incident with a fault-free $1$-dimension edge.

 (a) $(b^1, a^0)\in F$. Since $|F\cap \partial D_1|=2$ and $(b^1, a^0)\in F$, $b^0$ is incident with at least one fault-free $1$-dimension edge, say, $(b^0, a^3)$. Let $H_0=\langle b^0, a^0, v, u\rangle$ and $(a^1, b^2), (a^2, b^3)$ be fault-free $1$-dimension edges. Since $|F^2|=|F^3|=0$, there is a Hamiltonian path $H_i$ in $BH^{1, i}_1$ that joins $a^i$ and $b^i$ for $i=2, 3$.  Hence the cycle $C=\langle u, b^1, H_1, a^1, b^2, H_2, a^2, b^3, H_3, a^3, b^0, H_0, u\rangle$ (see figure \ref{appendix} (d)) is the desired cycle.

 (b) $(b^1, a^0)\notin F$.  Since $|F\cap D_1|=2$, $v$ or $b^0$ is incident with a fault-free $1$-dimension edge. Without loss of generality, we can assume that $b^0$ is incident with a fault-free $1$-dimension edge $(b^0, a^3)$. Let $H_0=\langle b^0, u, v, a^0\rangle$ and $(a^1, b^2), (a^2, b^3)$ be fault-free $1$-dimension edges. Since $|F^2|=|F^3|=0$, there is  a Hamiltonian path $H_i$ in $BH^{1, i}_{1}$ that joins $a^i$ and $b^i$ for $i= 2, 3$. Thus, the cycle $C=\langle b^0, H_0, a^0, b_1, H_1, a^1, b^2, H_2, a^2, b^3, H_3, a^3, b^0\rangle$ (see figure \ref{appendix} (e)) is the desired cycle.

 \noindent{\bf Subcase 2.1.3: } $|F^2|=1$.

 We can represent $BH^{1, 0}_1$  as $\langle u, v, a^0, b^0, u\rangle$. Since $|F\cap \partial D_1|=2$, we can see that there is an edge $(x, y)$ in $\{(v, a^0), (u, b^0)\}$ such that $x$ (resp.,  $y$) is incident with a  fault-free $j$-dimension edge. Without loss of generality, we can assume that $(x, y)=(u, b^0)$ and $(u, b^1), (b^0, a^3)$ are fault-free $1$-dimension edges.
 Let $H_1=\langle b^0, a^0, v, u\rangle$ and $(a^2, b^2)$ be the faulty edge in $BH^{1, 2}_1$. Since $\delta (BH_2-F)\ge 2$, $a^2$ (resp.,  $b^2$) is incident with a fault-free $1$-dimension edge $(a^2, b^3)$ (resp.,  $(a^1, b^2)$). Note that $|F^1|=|F^3|=0$. There is a Hamiltonian path $H_i$ in $BH^{1, i}_1$ that joins $a^i$ and $b^i$ for $i=1, 3$. Thus, the cycle $C=\langle u, b^1, H_1, a^1, b^2, H_2, a^2, b^3, H_3, a^3, b^0, H_0, u\rangle$ (see figure \ref{appendix} (f)) is the desired cycle.

\noindent{\bf Subcase 2.2: } $e\in \partial D_1$.

Without loss of generality, we can assume that $e=(u, v)=(a^0, b^1)$.

\noindent{\bf Subcase 2.2.1: } $|F^0|=1$(or $|F^1|=1$).

Suppose that $V(BH^{1, i}_1)=\{ a^i, b^i, c^i, d^i\}$.

\noindent{\bf Subcase 2.2.1.1: } $(a^0, x)$ is a faulty edge, where $x\in\{b^0, d^0\}$.

Without loss of generality, we can assume that $x=b^0$.
Since $\delta(BH_2-F)\ge 2$, $b^0$ is incident with a fault-free $1$-dimension edge, say, $(b^0, a^3)$. Let $H_0=BH^{1, 0}_1-(a^0, b^0)$ and $(a^1, b^2), (a^2, b^3)$ be fault-free $1$-dimension edges. Note that $|F^1|=|F^2|=|F^3|=0$. There is a Hamiltonian path $H_i$ of $BH^{1, i}_1$ that joins $a^i$ and $b^i$ for $i=1, 2, 3$. Thus, the cycle $C=\langle u, v, H_1, a^1, b^2, H_2, a^2, b^3, H_3, a^3, b^0, H_0, u\rangle$ (see figure \ref{appendix} (g)) is the desired cycle.

\noindent{\bf Subcase 2.1.1.2: } $(c^0, x)$ is a faulty edge where $x\in\{b^0, d^0\}$.

Without loss of generality, we can assume that $x=b_0$.

(a) Suppose that each vertex in $BH_2$ is incident with at least one fault-free $1$-dimension edge. Then,
there are two disjoint fault-free $1$-dimension edges between $BH^{1, i}_1$ \vspace{0.5ex}and $BH^{1, i+1}_1$ for all $i=1, 2, 3$, namely $(a^i, b^{i+1}), (c^i, d^{i+1})$ for $i=1, 2, 3$.
 Since $\delta(BH_2-F)\ge2$ and $(b^0, c^0)\in F$,  $(c^0, b^1)$ or $(c^0, d^1)$ is a fault-free $1$-dimension edge.

  If $(c^0, d ^1)\notin F$, the cycle $\langle u, v, c^0, d^0, c^3, d^3, c^2, d^2, c^1, d^1, a^1, b^2, a^2, b^3, a^3, b^0, u\rangle$ (see figure \ref{appendix} (h)) is the desired cycle.

  If $(c^0, d^1)\notin F$, the cycle $\langle u, v, c^1, d^2, c^2, d^3, c^3, d^0, c^0, d^1, a^1, b^2, a^2, b^3, a^3, b^0, u\rangle$ (see figure \ref{appendix} (i)) is the desired cycle.

(b) Suppose that there is a vertex $y$ such that $y$ is incident with two faulty $1$-dimension edges. Without loss of generality, we can let $y=a^1$. Then, the cycle $C=\langle u, v, a^1, d^1, c^0, d^0, c^3, d^3, c^2, d^2, c^1, b^2,$ $a^2, b^3, a^3, b^0, u\rangle$ (see figure \ref{appendix} (j)) is the desired cycle.

\noindent{\bf Subcase 2.2.2: } $|F^2|=1$(or $|F^3|=1$).

Let $(a^2, b^2)$ be the faulty edge in $BH^{1, 2}_1$. Since $\delta(BH_2-F)\ge 2$, $a^2$ (resp.,  $b^2$) is incident with a fault-free $1$-dimension edge $(a^2, b^3)$ (resp.,  $(a^1, b^2)$). Since $|F\cap\partial D_1|=2$, there is a fault-free $1$-dimension edge incident with $b^0$ or $d^0$. Without loss of generality, we can assume that $(b^0, a^3)$ is a fault-free $1$-dimension edge. Let $H_i=BH^{1, i}_1-(a^i, b^i)$ for $i=0, 1, 2, 3$. Then, the cycle $C=\langle u, v, H_1, a^1, b^2, H_2, a^2, b^3, H_3, a^3, b^0, H_0, u\rangle$ (see figure \ref{appendix} (k)) is the desired cycle.

\section{Conclusion}

In this paper, we consider the problem  proposed by  Zhou \cite{Qingguozhou} and show that each fault-free edge lies on a fault-free Hamiltonian cycle after no more than  $4n-5$ faulty edges occur when each vertex is incident with at least two fault-free edges. Our result is optimal because there is a counterexample  when $|F|=4n-4$ as follows. Let $u=(0, 0, \cdots, 0), v=(2, 0, \cdots, 0), x=(1, 0, \cdots, 0), y=(3, 0, \cdots, 0)$ and $F=\{(u, a)\mid a\in N_{BH_n}(u)\setminus\{x, y\}\}\bigcup \{(v, b)\mid b\in N_{BH_n}(u)\setminus\{x, y\}\}$. Then, $|F|=4n-4$. Suppose that there is a Hamiltonian cycle $C$. Note that $u, v\in V(C)$.  Then, $C$  contains $(u, x), (u, y), (v, x), (v, y)$. Note that $\langle u, x, v, y, u\rangle$ is a cycle of length $4$, which is a contradiction. Thus, there is no Hamiltonian cycle in $BH_n-F$.

\section*{Acknowledgement}
This research is supported by the National Natural Science Foundation of China (11571044, 61373021),  the Fundamental Research Funds for the Central University of China.

\end{CJK*}

\end{document}